\documentclass[a4paper,12pt]{amsart}
\usepackage{wrapfig}
\usepackage[dvips]{graphicx}
\usepackage{amsmath,amsthm,amssymb}

\theoremstyle{definition}
\newtheorem{thm}{Theorem}[section]
\theoremstyle{definition}

\theoremstyle{definition}
\newtheorem{lem}[thm]{Lemma}
\theoremstyle{definition}

\theoremstyle{definition}
\newtheorem{prop}[thm]{Proposition}
\theoremstyle{definition}

\theoremstyle{definition}

\theoremstyle{definition}

\theoremstyle{definition}

\usepackage{amscd}

\makeindex
\title{Circle correspondence $C^*$-algebras}
\author{Shinji Yamashita}
\address{Graduate School of Mathmatics, Kyushu University,
Hakozaki, Fukuoka, 812-8581, JAPAN}
\email{s-yamashita@math.kyushu-u.ac.jp}
\date{}

\begin{document}
\maketitle

\begin{abstract}
We investigate Cuntz-Pimsner $C^*$-algebras associated with certain correspondences of the unit circle $\mathbb{T}$. We analyze these $C^*$-algebras by analogy with irrational rotation algebras $A_\theta$ and Cuntz algebras $\mathcal{O}_n$. We construct a Rieffel type projection, study the fixed point algebras of certain actions of finite groups, and calculate the entropy of a certain endomorphism. We also study the induced map of the dual action of the gauge action on $K$-groups.
\end{abstract}


\section{Introduction}
\label{sec:Introduction}

In \cite{KajiwaraWatatani}, Kajiwara and Watatani introduced a $C^*$-algebra $\mathcal{O}_R $ associated with a rational function $R$ as a Cuntz-Pimsner algebra of a Hilbert bimodule over $A=C(J_R)$, which is the algebra of continuous functions on Julia set $J_R$ of $R$. They proved that if the degree of $R$ is at least two, then a $C^*$-algebra $\mathcal{O}_R$ is a \textit{Kirchberg algebra} (purely infinite, simple, nuclear, separable $C^*$-algebra) satisfying UCT. In this framework, one of the most important issues is to observe relations between $C^*$-algebras $\mathcal{O}_R $ and complex dynamical systems. But it is also important to examine the properties for elementary rational functions. \\
The Cuntz algebras and the irrational rotation algebras have been examined by many authors. These algebras are simple and universal $C^*$-algebras with certain commutation relations. They have their own properties, which are not easily observed in more general $C^*$-algebras.

In this paper, we deal with a $C^*$-algebra $\mathcal{O}_{(m,n)}(\mathbb{T})$ which is generated by two elements $z$ and $S_1$ with commutation relations $z^nS_1=S_1z^m, \sum_{i=0}^{n-1} z^iS_1S_1^*z^{-i}=1$. This $C^*$-algebra include the above $C^*$-algebra $\mathcal{O}_R$ associated with the elementary rational function $R(z)=z^n$ on $J_R=\mathbb{T}$. Our algebras have appeared in several papers (for example, a groupoid $C^*$-algebra \cite{Deaconu} and a topological graph $C^*$-algebra \cite{Katsura}). The main purpose of the present paper is to examine specific properties of our algebras like Cuntz algebras $\mathcal{O}_n$ and irrational rotation algebras $A_\theta$. In fact, our algebras contain Cuntz algebras as subalgebras and are generated by two elements like the $C^*$-algebra $A_\theta$.

This paper is organized as follows. In section \ref{sec:Preliminaries}, we give some preliminaries. In section \ref{sec:Projection}, we construct a projection on matrix algebra over $\mathcal{O}_{(m,n)}(\mathbb{T}) $ for $(m,n)=(1,2)$. This construction is similar to that of the $C^*$-algebra $A_\theta$ presented by Rieffel \cite{Rieffel}. In section \ref{sec:Subalgebra}, we discuss the $C^*$-subalgebras of $\mathcal{O}_{(m,n)}(\mathbb{T}) $. In the case of the $C^*$-algebra $A_\theta$, for any positive integer $k$, $A_{k\theta}$ is naturally a $C^*$-subalgebra of $A_\theta$. We treat an analog of this problem for $\mathcal{O}_{(m,n)}(\mathbb{T}) $. We apply this result to obtain the fixed point algebra of a cyclic group action on $\mathcal{O}_{(m,n)}(\mathbb{T}) $. Moreover we determine the fixed point of a symmetry action. In section \ref{sec:Entropy}, we calculate the entropy of a certain endomorphism that seems like Cuntz's canonical endomorphism. In section \ref{sec:Dual}, we consider the gauge action and its dual action of $\mathcal{O}_{(m,n)}(\mathbb{T}) $. In particular, we investigate the induced map of the dual action of the gauge action on $K$-groups. As a corollary, we obtain another computation of $K$-groups. 


\section{Preliminaries}
\label{sec:Preliminaries}
In this section, we recall the Cuntz-Pimsner algebras and introduce the $C^*$-algebras $\mathcal{O}_{(m,n)}(\mathbb{T})$. A \textit{(right) Hilbert $A$-module} $X$ is a Banach space (whose norm is $\| \cdot \|_2$) with a right action of a $C^*$-algebra $A$ and an $A$-valued inner product $\langle \cdot , \cdot \rangle $ satisfying:
\begin{enumerate}
\item{ $\langle f, ga \rangle =\langle f, g \rangle a$}
\item{ $\langle f, g \rangle^* =\langle g,f \rangle $}
\item{ $ \langle f,f \rangle \ge 0$ and $\| f \|_2=\| \langle f,f \rangle \|^{1/2}$}
\end{enumerate}
for $f,g \in X$ and $a \in A$. The Hilbert $A$-module is \textit{finitely generated} if there exists $\{ u_i \}_{i=1}^n \subset X$ such that $f=\sum_{i=1}^n u_i \langle u_i ,f \rangle$ for any $f \in X$. A Hilbert $A$-module $X$ is \textit{full} if $\overline{\mathrm{span}}\{ \langle f,g \rangle | f,g \in X \} =A$. We only use finitely generated full Hilbert bimodules. 
For a Hilbert $A$-module $X$, we denote by $L(X)$ the $C^*$-algebra of all adjointable operators on $X$. For $f,g \in X$, the rank-one operator $\theta_{f,g} \in L(X)$ is defined by $\theta_{f,g}(h)=f\langle g,h \rangle$ for $h \in X$. Set $K(X)=\overline{\mathrm{span}}\{\theta_{f,g} | f,g \in X \}$. If $X$ is finitely generated by $\{ u_i\}_{i=1}^n$ and $\langle u_i,u_j \rangle = \delta_{i,j}1$, then $L(X)=K(X)$ and $K(X)$ is isomorphic to $M_n(A)$ via $\{u_i \}_{i=1}^n$. \\
We recall the Cuntz-Pimsner algebras. Let $A$ be a $C^*$-algebra and $X$ be a full Hilbert $A$-module that is finitely generated by $\{ u_i\}_{i=1}^n$ and let $\phi : A \longrightarrow L(X)$ be a faithful *-homomorphism. We shall define a left action by $a \cdot f:=\phi (a) f$ for $a \in A, f\in X$. Then the \textit{Cuntz-Pimsner algebra} $\mathcal{O}_X$ is the universal $C^*$-algebra generated by $A$ and $\{ S_f|f \in X \}$ with the following relation: for $a,b \in A, f,g \in X ,\alpha ,\beta \in \mathbb{C}$,
\begin{equation*}
S_{\alpha f+\beta g}=\alpha S_f + \beta S_g,\quad a\cdot S_f \cdot b=S_{a \cdot f \cdot b} ,\quad S_f^*S_g=\langle f,g \rangle ,\quad \sum_{i=1}^n S_{u_i}S_{u_i}^*=1.
\end{equation*}


Let us define the $C^*$-algebra $\mathcal{O}_{(m,n)}(\mathbb{T})$ for $m,n \in \mathbb{N}$ using the Cuntz-Pimsner construction. Set $\Omega_{(m,n)}=\{ (t^m,t^n) \in \mathbb{T} \times \mathbb{T} |t \in \mathbb{T} \}$. Set $X_{(m,n)}=C(\Omega_{(m,n)} )$. Then $X_{(m,n)}$ is a $C(\mathbb{T})$-$C(\mathbb{T})$ bimodule by 
\begin{equation*}
(a\cdot f \cdot b)(x,y)=a(x)f(x,y)b(y)
\end{equation*} 
for $a,b \in C(\mathbb{T}) ,f \in X_{(m,n)}, (x,y) \in \Omega_{(m,n)}$. We introduce a $C(\mathbb{T})$-valued inner product $\langle \cdot ,\cdot \rangle$ on $X_{(m,n)}$ by
\begin{equation*}
\langle f,g \rangle(y)=\sum_{\{ x \in \mathbb{T} | (x,y) \in \Omega_{(m,n)} \}} \overline{f(x,y)}g(x,y)
\end{equation*}
for $f,g \in X_{(m,n)}$ and $y \in \mathbb{T}$. Put $\| f\|_2=\| \langle f,f \rangle \|_\infty^{1/2}$. Then $X_{(m,n)}$ is a full Hilbert bimodule over $C(\mathbb{T})$ without completion (see Corollary 2.3 of \cite{KajiwaraWatatani}).
Let us denote the greatest common divisor of $m,n \in \mathbb{N}$ by $\mathrm{gcd}(m,n)$. In this case, the Hilbert module $X_{(m,n)}$ is a finitely generated by $u_i(x,y)=\frac{1}{\sqrt{n_0}}x^{i-1}\ (i=1,\cdots , n_0:=n/\mathrm{gcd}(m,n))$ and $\{ u_i \}_{i=1}^{n_0}$ satisfies $\langle u_i,u_j \rangle =\delta_{i,j}1$. Moreover $X_{(m,n)}$ is isomorphic to $X_{(m_0,n_0)}$ as a Hilbert module where $m_0=m/\mathrm{gcd}(m,n)$. \textbf{Thus we always assume $\mathbf{\mathrm{\textbf{gcd}}(m,n)=1}$.}

We denote $S_i=S_{u_i}$ for simplicity. Let $z$ be an element of $C(\mathbb{T})$ defined by $z(x)=x$ for any $x \in \mathbb{T}$. Then $\mathcal{O}_{(m,n)}(\mathbb{T})$ is the universal $C^*$-algebra generated by a (full-spectrum) unitary $z$ and isometries $S_1,\cdots ,S_n $ with the relation
\begin{equation*}
zS_i=S_{i+1}\ (1 \le i \le n-1), \quad zS_n=S_1z^m, \quad \sum_{i=1}^nS_iS_i^*=1.
\end{equation*}
We will often use this relation. From this relation, we notice that $\mathcal{O}_{(m,n)}(\mathbb{T})$ is generated by the two elements $z$ and $S_1$.\\
For $k \in \mathbb{N}_0:=\mathbb{N} \cup \{ 0 \}$, we define the set $\mathcal{W}_n^{(k)}$ of $k$-triples by $\mathcal{W}_n^{(0)} = \{ \emptyset \}$ and $\mathcal{W}_n^{(k)}=\{ (i_1,i_2,\cdots , i_k)| i_j \in \{ 1,\cdots ,n \} \}$. Set $\mathcal{W}_n=\bigcup_{k=0}^\infty \mathcal{W}_n^{(k)}$. For $\mu =(i_1,\cdots ,i_k ) \in \mathcal{W}_n$, we denote its length $k$ by $| \mu |$ and set $S_\mu = S_{i_1}S_{i_2} \cdots S_{i_k}$. Note that $|\emptyset |=0, S_{ \emptyset }=1$. For $\mu = (i_1,\cdots , i_k), \nu=(j_1,\cdots , j_l) \in \mathcal{W}_n$, we define their product $\mu \nu \in \mathcal{W}_n$ by $\mu \nu =(i_1,\cdots i_k , j_1 , \cdots j_l )$. 
Using these notations, since $X_{(m,n)}$ is finitely generated, $\mathcal{O}_{(m,n)}(\mathbb{T})$ is presented by
\begin{equation*}
\mathcal{O}_{(m,n)}(\mathbb{T})= \overline{\mathrm{span}}\Bigl\{ S_\mu z^k S_\nu^* \Bigl| k \in \mathbb{Z} ,\mu , \nu \in \mathcal{W}_n \Bigr\} .
\end{equation*}

There exists an action $\alpha : \mathbb{T} \ni t \longmapsto \alpha_t \in \mathrm{Aut}(\mathcal{O}_{(m,n)}(\mathbb{T}))$ with $\alpha_t(S_f)=tS_f$ that is called the \textit{gauge action}. The fixed point algebra $\mathcal{O}_{(m,n)}(\mathbb{T})^\alpha$ by the gauge action is the \textit{$(m,n)$-type Bunce-Deddens algebra} $\mathcal{B}_{(m,n)}:= \varinjlim_{k \ge 0} \{ M_{n^k}(C(\mathbb{T} )) , \phi_k \}$ where $\phi_k : M_{n^k}(C(\mathbb{T} )) \longrightarrow M_{n^{k+1}}(C(\mathbb{T} ))$ is defined by
\[
\phi_k(a)=
\left(
\begin{array}{cccc|c}
a&\cdots&0\\
\vdots&\ddots^{\times n}&\vdots\\
0&\cdots&a
\end{array}
\right)
\ (a \in M_n(\mathbb{C} )), \quad
\phi_k
\left(
\begin{array}{cccc|c}
0 & z \\
1_{n^k-1} & 0
\end{array}
\right)
=
\left(
\begin{array}{cccc|c}
0 & z^m \\
1_{n^{k+1}-1} & 0
\end{array}
\right).
\]
The element $S_\mu z^k S_\nu^* \in \mathcal{O}_{(m,n)}(\mathbb{T})^\alpha$ corresponds to $z^k \otimes e_{\mu \nu}=z^k \otimes e_{\mu_1 \nu_1} \otimes \cdots \otimes e_{\mu_{|\mu |} \nu_{|\nu |}}$ where $e_{ij}$ is the $(i,j)$-matrix unit.\\
We shall discuss about some properties on $\mathcal{O}_{(m,n)}(\mathbb{T})$. Recently, Katsura showed that some his algebras, so-called topological graph $C^*$-algebras, are Kirchberg algebras satisfying UCT. Our algebras $\mathcal{O}_{(m,n)}(\mathbb{T})$ are contained in his algebras. Katsura also computed these $K$-groups in Appendix A of \cite{Katsura3}.


\begin{thm}[Katsura\ \cite{Katsura3}]
\label{thm:Katsura}\normalfont\slshape
Suppose that $\mathrm{gcd}(m,n)=1$.
\begin{enumerate}
\item{
For $m \ge 1,n \ge 2$, $\mathcal{O}_{(m,n)}(\mathbb{T})$ are Kirchberg algebras satisfying UCT.
}
\item{
\begin{enumerate}
\item{For $n \ge 2$, $K_0(\mathcal{O}_{(1,n)}(\mathbb{T} ))=\mathbb{Z} \oplus \mathbb{Z}_{n-1}$, $K_1(\mathcal{O}_{(1,n)}(\mathbb{T} ))=\mathbb{Z}$.}
\item{For $m \ge 2$, $K_0(\mathcal{O}_{(m,1)}(\mathbb{T} ))=\mathbb{Z}$, $K_1(\mathcal{O}_{(m,1)}(\mathbb{T} ))=\mathbb{Z} \oplus \mathbb{Z}_{m-1}$}.
\item{For $m,n \ge 2$, $K_0(\mathcal{O}_{(m,n)}(\mathbb{T}))=\mathbb{Z}_{n-1}$, $K_1(\mathcal{O}_{(m,n)}(\mathbb{T}))=\mathbb{Z}_{m-1}$}
\end{enumerate}
}
\end{enumerate}
where $\mathbb{Z}_k$ is the cyclic group $\mathbb{Z}/k\mathbb{Z}$.
\end{thm}

\begin{proof}
Since we use the computation of $K$-groups in Section \ref{sec:fixsym}, we give a brief discussion. For the general Cuntz-Pimsner algebras $\mathcal{O}_X$ associated with $(A,X,\phi )$, Katsura gave the six-term exact sequence (without the KK-theoretic method):
\[
\begin{CD}
K_0(A) @>\mathrm{id}-[X]_0>> K_0(A) @>\iota_*>> K_0(\mathcal{O}_X )\\
@A\delta_1AA @. @VV\delta_0V \\
K_1(\mathcal{O}_X) @<\iota_*<< K_1(A) @<\mathrm{id}_*-[X]_1<< K_1(A)
\end{CD}
\]
where $\iota : A \longrightarrow \mathcal{O}_X$ is the natural inclusion and for the case where $X$ is finitely generated by $\{ u_i \}_{i=1}^n$ satisfying $\langle u_i , u_j \rangle =\delta_{i,j}$, the map $[X]_i$ are the following composition maps:
\begin{equation*}
[X]_i: \quad K_i(A) \overset{\phi_*}{ \longrightarrow } K_i(M_n(A)) \overset{\cong}{ \longrightarrow } K_i(A).
\end{equation*}
Now, we consider our algebras $\mathcal{O}_{(m,n)}$. Since $K_0(C(\mathbb{T} ))=\mathbb{Z}[1]_0, K_1(C(\mathbb{T} ))=\mathbb{Z}[z]_1$, and for the left action $\phi : C(\mathbb{T} ) \longrightarrow L(X_{(m,n)}) \cong M_n(A)$,
\[
\phi (1)=
\left(
\begin{array}{cccc|c}
1&\cdots&0\\
\vdots&\ddots^{\times n}&\vdots\\
0&\cdots&1
\end{array}
\right)
,\quad
\phi (z)
=
\left(
\begin{array}{cccc|c}
0 & z^m \\
1_{n-1} & 0
\end{array}
\right).
\]
Hence, $[\phi (1)]_{0,M_n(C(\mathbb{T} ))}=n[1]_{0,C(\mathbb{T} ))}$,\ $[\phi (z)]_{1,M_n(C(\mathbb{T} ))}=m[z]_{1,C(\mathbb{T} ))}$. Consequently, $[X_{(m,n)}]_0$ is an $n$-times map and $[X_{(m,n)}]_1$ is an $m$-times map on $\mathbb{Z}$.
\end{proof}

When $m \ge 2 ,n=1$ for $\mathcal{O}_{(m,n)}(\mathbb{T})$, the $C^*$-algebras $\mathcal{O}_{(m,1)}(\mathbb{T})$ are not Kirchberg algebras: in fact, they are not simple. Moreover, these algebras are transformation group $C^*$-algebras on solenoid groups and they have been systematically examined by Brenken-J\o rgensen \cite{BJ}, Brenken \cite{Bre}. We recall the solenoid group. Define
\begin{equation*}
S_m=\Bigl\{ (x_i)_{i=0}^\infty \in \prod_{i=0}^\infty \mathbb{T} \ \Bigl| x_{i+1}^m = x_i (i \in \mathbb{N}_0) \Bigr\} .
\end{equation*}
Then $S_m$ be a compact connected abelian group: it is called the \textit{solenoid group}. We refer to \cite{Williams2} for the solenoid. Let us define a group automorphism $\sigma$ on $S_m$ by $\sigma (x)_i=x_{i+1}$ for $x \in S_m$. Then $\mathcal{O}_{(m,1)}(\mathbb{T} )$ is isomorphic to the crossed product $C(S_m) \rtimes_\sigma \mathbb{Z}$ \cite{KajiwaraWatatani2}. For completeness, we give a brief proof for some properties on solenoid $C^*$-algebras examined by Brenken-J\o rgensen \cite{BJ}, Brenken \cite{Bre}.

\begin{thm}[Brenken-J\o rgensen \cite{BJ}, Brenken \cite{Bre}]
\normalfont\slshape
\ \\
\ For $m \ge 2$, the solenoid $C^*$-algebra $\mathcal{O}_{(m,1)}(\mathbb{T}) \cong C(S_m) \rtimes_\sigma \mathbb{Z}$ is NGCR, AF-embeddable, non-simple, residually finite dimensional.
\end{thm}

\begin{proof}
For $x \in S_m$, define $O(x)=\{ \sigma^k(x) \in S_m | k \in \mathbb{Z} \}$ and $D_m=\{ x \in X | O(x) \mbox{ is dense in } S_m  \}$. Then $D_m$ is a dense set of $S_m$. Hence, there exist $x, y \in S_m$ such that $O(x) \ne O(y)$ and $\overline{O(x)} = S_m = \overline{O(y)}$, where $\overline{O(x)}$ is the closure set of $O(x)$ in $S_m$ and the consequence of this observation induces the orbit space of the dynamics $(S_m , \sigma )$ to not be a $T_0$-topological space. Hence, $C(S_m) \rtimes_\sigma \mathbb{Z}$ is NGCR (see Section 8 of \cite{Williams}).\\
For $k \ge 1$, let $\mathrm{ Per}_k(\sigma )=\{ x \in S_m| \sigma^k(x)=x , \sigma^l(x) \ne x (1 \le l < k )\}$ be the set of $k$-period points. Then for any $k \ge 1$, $\mathrm{Per}_k(\sigma )$ is not empty, and moreover $\mathrm{Per}(\sigma ):= \cup_{k=1}^\infty \mathrm{Per}_k(\sigma )$ is a countable dense set in $S_m$. Hence the non-wandering set of $(S_m, \sigma )$ coincides with $S_m$. These results imply that $C(S_m) \rtimes_\sigma \mathbb{Z}$ is AF-embeddable according to the a work of Pimsner \cite{Pimsner2} and non-simple (the existence of a periodic points implies that $\sigma$ is not minimal).\\
Let us state that $C(S_m) \rtimes_\sigma \mathbb{Z}$ is residually finite dimensional. We shall show that there is a countable family $\{ \pi_n \}$ of the representations for finite dimension $C^*$-algebras such that $\pi := \oplus_n \pi_n$ is faithful. For each $x \in \mathrm{Per}(\sigma )$, let $k_x$ be the period of $x$ and let us define $\rho_x : C(S_m) \longrightarrow M_{k_x}(\mathbb{C} )$ by $\rho_x(f)=\mathrm{diag}(f(x), f(\sigma (x)), \cdots , f(\sigma^{{k_x}-1}(x)))$ and for each $z \in \mathbb{T}$, let us define a unitary $u_{x,z}$ by 
\[u_{x,z}=
\left(
\begin{array}{cccc|c}
0 & z \\
1_{k_x-1} & 0
\end{array}
\right).
\]
Then $(\rho_x , u_{x,z} )$ satisfies the covariance relation, and we denote the covariance representation by $\pi_{x,z} := \rho_x \rtimes u_{x,z} : C(S_m) \rtimes_\sigma \mathbb{Z} \longrightarrow M_{k_x}(\mathbb{C} )$. Set
\begin{equation*}
\pi_z := \bigoplus_{x \in \mathrm{Per}(\sigma )} \pi_{x,z} : C(S_m) \rtimes_\sigma \mathbb{Z} \longrightarrow \bigoplus_{x \in \mathrm{Per}(\sigma )} M_{k_x}(\mathbb{C} ).
\end{equation*}
Let $\{ z_l \}_{l=1}^\infty$ be a dense set of $\mathbb{T}$. Set $\pi:=\oplus_{l=1}^\infty \pi_l$. We shall show that $\pi$ is a faithful representation.\\
Let $a \in C(S_m) \rtimes_\sigma \mathbb{Z}$ be a positive element such that $\pi (a)=0$: this implies $\pi_{z_l}(a)=0$ for $l \in \mathbb{N}$. Since $\mathbb{T} \ni z \longmapsto \pi_z(a)$ is continuous, $\pi_z(a)=0$ for any $z \in \mathbb{T}$. For $w \in \mathbb{T}$, let us define an automorphism $ \lambda_w^{(x)}$ on $M_{k_x}(\mathbb{C} )$ by $\lambda_w^{(x)}(e_{ij})=w^{i-j}e_{ij}$, where $e_{ij}$ is the $(i,j)$-matrix unit. Set $\lambda_w:= \oplus_{x \in \mathrm{Per}(\sigma )} \lambda_w^{(x)}$. Let $\Psi : C(S_m) \rtimes_\sigma \mathbb{Z} \longrightarrow C(S_m)$ be the (canonical) faithful conditional expectation. Then
\begin{equation*}
0=\int_{\mathbb{T}}\lambda_z \bigl(\int_\mathbb{T} \lambda_w (\pi_z(a))dw \bigr) dz = \bigoplus_{x \in \mathrm{Per}(\sigma )} \rho_x(\Psi(a)).
\end{equation*}
Since $\mathrm{Per}(\sigma )$ is dense in $S_m$, $\Psi(a)=0$ and since $\Psi$ is faithful, we can conclude that $a=0$. Consequently, $\pi$ is a faithful representation.
\end{proof}


\section{Rieffel-type projection on matrix algebra over $\mathcal{O}_{(1,2)}(\mathbb{T}) $}
\label{sec:Projection}

In \cite{Rieffel}, Rieffel explicitly described some projections in the irrational rotation algebras $A_\theta$ and obtained the value that the trace has on them.
In this section, using a similar method to \cite{Rieffel}, we would like to construct a projection on $M_2(\mathcal{O}_{(1,2)}(\mathbb{T}) )$ which is not von Neumann equivalent to 1.

\subsection{Construction of projection}

We shall construct a projection of $\mathcal{O}_{(1,2)}(\mathbb{T})$. Define a *-homomorphism $\phi :C(\mathbb{T}) \longrightarrow C(\mathbb{T})$ by $\phi (a)(t)=a(2t)$ for $a \in C(\mathbb{T}), t \in \mathbb{T}=[0,1)$. Let us define an element on $M_2(\mathcal{O}_{(1,2)}(\mathbb{T}) )$ by
\begin{eqnarray*}
P=
\begin{pmatrix}
P_{11} & P_{12}\\
P_{21} & P_{22} \\
\end{pmatrix}
=\begin{pmatrix}
S_1a_1 + \phi (a_0) +a_1S_1^* &S_2a_1 + b_1S_2^* \\
S_2b_1 + a_1S_2^* & S_1b_1 + \phi (b_0) +b_1S_1^*, \\
\end{pmatrix}
\end{eqnarray*}
where $a_0,a_1,b_0,b_1 \in A$ are real-valued functions. Then $P$ is a self-adjoint from the construction. We will show the following theorem in this section.

\begin{thm}\normalfont\slshape
\label{thm:projection}
Suppose that $(m,n)=(1,2)$. Then we can construct a projection $P$ of $M_2(\mathcal{O}_{(1,2)}(\mathbb{T})  )$ that is not von Neumann equivalent to 1 and $[P]_0$ is $-4$ of the group $K_0(\mathcal{O}_{(1,2)}(\mathbb{T}) )=\mathbb{Z}$. Let $\varphi$ be the $\log 2$-KMS state for the gauge action (see \cite{KajiwaraWatatani}) and $\tau$ be a normalized trace on $M_2(\mathbb{C} )$. Then $\varphi \otimes \tau (P)=\frac{7}{16}$.
\end{thm}

\begin{proof}
We need to investigate how to impose conditions in order to satisfy $P^2=P$.\\
First, we shall see that $P_{11}=P_{11}P_{11}+P_{12}P_{21}$. The right-side term becomes
\begin{eqnarray*}
\lefteqn{P_{11}P_{11}+P_{12}P_{21}}\hspace{1cm}\\
&=&\phi (a_0^2)+a_1^2 +b_1^2 +\phi(a_1^2) +
S_1\bigl( a_1(a_0+\phi (a_0)) \bigr) + \bigl( a_1(a_0+\phi (a_0)) \bigr)S_1^* \\
& &+\phi (a_1) \phi^2 (a_1)S_1S_1 +\phi (a_1)\phi^2(b_1) S_2S_2+ S_1^*S_1^* \phi (a_1) \phi^2 (a_1)+S_2^*S_2^* \phi (a_1)\phi^2(b_1)\\
\end{eqnarray*}
To obtain the equation $P_{11}=P_{11}P_{11}+P_{12}P_{21}$, we have to impose the following conditions:
\begin{enumerate}
 \item{$\phi (a_0) -\phi (a_0^2)=a_1^2+b_1^2+\phi (a_1^2)$}
 \item{$a_1(a_0+\phi (a_0))=a_1$}
 \item{$a_1\phi (a_1)=0,a_1 \phi (b_1)=0$}. 
\end{enumerate}
Let us take $a_0$ to be
\[
a_0(t) =
\left\{
  \begin{array}{cl}
  4\bigl( t-\frac{1}{2} \bigr) & \mbox{if $t \in \bigl[ \frac{1}{2}, \frac{3}{4} \bigr]$} \\
  1-8\bigl( t-\frac{3}{4} \bigr) & \mbox{if $t \in \bigl[ \frac{3}{4}, \frac{7}{8} \bigr]$} \\
  0 & \mbox{otherwise}
  \end{array}
\right.
\]
and also define $a_1,b_1$ by
\[
a_1(t) =
\left\{
  \begin{array}{cl}
  \sqrt{\phi (a_0)(t)-\phi (a_0^2)(t)} & \mbox{if $t \in \bigl[ \frac{3}{4}, \frac{7}{8} \bigr]$} \\
  0 & \mbox{otherwise}
  \end{array}
\right.
b_1(t) =
\left\{
  \begin{array}{cl}
  \sqrt{\phi (a_0)(t)-\phi (a_0^2)(t)} & \mbox{if $t \in \bigl[ \frac{1}{4}, \frac{3}{8} \bigr]$} \\
  0 & \mbox{otherwise}
  \end{array}
\right.
\]
\begin{figure}
 \centering
 \includegraphics[width=15cm,clip]{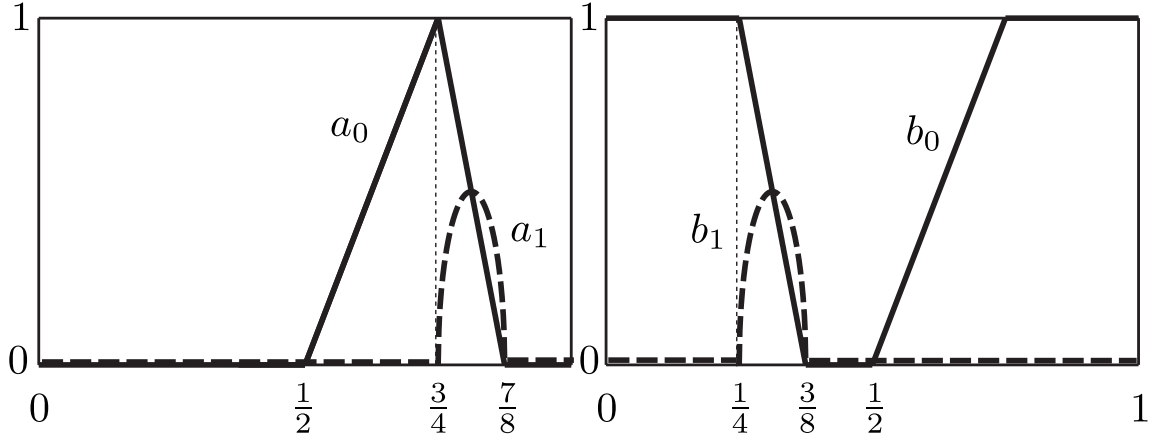}
 \caption{Graphs of functions $a_0,a_1,b_0,b_1$.}
 \label{fig:graph}
\end{figure}
We can check that these functions satisfy all the conditions.\\
Next, we shall check that $P_{22}=P_{21}P_{12}+P_{22}P_{22}$. 
\begin{eqnarray*}
\lefteqn{P_{21}P_{12}+P_{22}P_{22}}\hspace{1cm}\\
&=&\phi (b_0^2)+a_1^2 +b_1^2 +\phi(b_1^2) +
S_1\bigl( b_1(b_0+\phi (b_0)) \bigr) + \bigl( b_1(b_0+\phi (b_0)) \bigr) S_1^* \\
& &+\phi (b_1) \phi^2 (b_1)S_1S_1 +\phi (b_1)\phi^2(a_1) S_2S_2+ S_1^*S_1^* \phi (b_1) \phi^2 (b_1)+S_2^*S_2^* \phi (b_1)\phi^2(a_1)\\
\end{eqnarray*}
To obtain the equation $P_{22}=P_{21}P_{12}+P_{22}P_{22}$, we have to impose the following conditions:
\begin{enumerate}
 \item{$\phi (b_0) -\phi (b_0^2)=a_1^2+b_1^2+\phi (b_1^2)$}
 \item{$b_1(b_0+\phi (b_0))=b_1$}
 \item{$b_1\phi (b_1)=0,b_1 \phi (a_1)=0$}. 
\end{enumerate}
To satisfy these conditions, we define $b_0$ by
\[
b_0(t) =
\left\{
  \begin{array}{cl}
  1 & \mbox{if $t \in \bigl[ 0, \frac{1}{4} \bigr]$} \\
  1-8\bigl( t-\frac{1}{4} \bigr) & \mbox{if $t \in \bigl[ \frac{1}{4}, \frac{1}{2} \bigr]$} \\
  0 & \mbox{if $t \in \bigl[ \frac{1}{4}, \frac{1}{2} \bigr]$} \\
  4\bigl( t-\frac{1}{2} \bigr) & \mbox{if $t \in \bigl[ \frac{1}{2}, \frac{3}{4} \bigr]$} \\
  1 & \mbox{if $t \in \bigl[  \frac{3}{4}, 1 \bigr]$} \\
  \end{array}
\right.
\]
Next, we have to check that the off-diagonal part $P_{12}=P_{11}P_{12}+P_{12}P_{22},P_{21}=P_{21}P_{11}+P_{22}P_{21}$ is also affirmative. 
\begin{equation*}
P_{11}P_{12}+P_{12}P_{22}=S_2a_1\bigl( a_0+\phi (b_0) \bigr) + \bigl( \phi(a_0 )+b_0 \bigr) b_1S_2^*
\end{equation*}
and to compare the both terms, we need to impose $a_0+\phi (b_0)=1$ on $\mathrm{supp}(a_1) = [ 3/4,7/8]$ and $\phi (a_0)+ b_0=1$ on $\mathrm{supp}(b_1) = [ 1/4,3/8]$. But this equation is correct from the definition of functions. Consequently we conclude that $P$ is a projection on $M_2(\mathcal{O}_{(1,2)}(\mathbb{T}) )$.

\subsection{Proof of $[P]_0=-4$ in $K_0(\mathcal{O}_{(1,2)}(\mathbb{T}) $) }
Let us examine $[P]_0=-4$ in $K_0(\mathcal{O}_{(1,2)}(\mathbb{T}) )=\mathbb{Z}$. We recall an exact sequence
\begin{equation*}
0 \longrightarrow K(F(X_{(1,2)})) \longrightarrow \mathcal{T}_{(1,2)}(\mathbb{T}) \overset{\pi }{\longrightarrow} \mathcal{O}_{(1,2)}(\mathbb{T})  \longrightarrow 0
\end{equation*}
where $ K(F(X_{(1,2)}))$ is the $C^*$-algebra generated by one-rank operators on Fock space $F(X_{(1,2)})$ and $\mathcal{T}_{(1,2)}(\mathbb{T})$ is the Toeplitz $C^*$-algebra (see \cite{Pimsner}). From the six-term exact sequence in the proof in Theorem \ref{thm:Katsura}, we can show that the exponential map $\delta_0 : K_0(\mathcal{O}_{(1,2)}(\mathbb{T}) ) \longrightarrow K_1(K(F(X_{(1,2)})) )$ is a group isomorphism. Let us observe $\delta_0([P]_0)$.

There exists a projection $Q \in L(F(X_{(1,2)}))$ such that $T_{u_1}T_{u_1}^*+T_{u_2}T_{u_2}^*+Q=1$. Set $T_i =T_{u_i} (i=1,2)$ and set an element $H$ in $M_2(\mathcal{T}_{(1,2)}(\mathbb{T}) )$ by
\begin{eqnarray*}
H=
\begin{pmatrix}
H_{11} & H_{12}\\
H_{21} & H_{22} \\
\end{pmatrix}
=\begin{pmatrix}
T_1a_1 + \phi (a_0) +a_1T_1^* &T_2a_1 + b_1T_2^* \\
T_2b_1 + a_1T_2^* & T_1b_1 + \phi (b_0) +b_1T_1^* \\
\end{pmatrix}.
\end{eqnarray*}
$H$ satisfies $\pi^{(2)}(H)=P$. From the definition of $\delta_0$, $\delta_0[P]_0=[\exp (2\pi i H)]_1$ in $K_1(K(F(X_{(1,2)})))$. We can check $H^2=H-\mathrm{diag}(\phi (a_1^2)Q, \phi (b_1^2)Q ).$ Moreover, we can calculate
\begin{equation*}
H^m=H-\mathrm{diag}\Bigl( \sum_{k=0}^{m-2} \phi (a_0^k) \phi (a_1^2)Q , \ \sum_{k=0}^{m-2} \phi (b_0^k) \phi (b_1^2)Q \Bigr)
\end{equation*}
for $m \ge 2$. Let us define $\Delta_1 =\chi_{\mathrm{supp}a_1}=\chi_{[\frac{3}{4} , \frac{7}{8} ]}, \Delta_2 =\chi_{\mathrm{supp}b_1}=\chi_{[\frac{1}{4} , \frac{3}{8} ]}$, where $\chi_S$ is a characteristic function on $S \subset \mathbb{T}$. Then $a_1^2=(\phi (a_0) -\phi (a_0^2))\Delta_1 =(a_0 -a_0^2)\Delta_1$ and 
\begin{equation*}
\sum_{k=0}^{m-2} a_0^k a_1^2 =\sum_{k=0}^{m-2} a_0^k (a_0-a_0^2)\Delta_1=(a_0-a_0^m)\Delta_1.
\end{equation*}
In the same way, we can show that
\begin{equation*}
\sum_{k=0}^{m-2} a_0^k a_1^2 =(b_0-b_0^m)\Delta_2.
\end{equation*}
Hence,
\begin{equation*}
H^m = H -\mathrm{diag} \Bigl( \phi \bigl((a_0-a_0^m)\Delta_1\bigr)Q ,\ \phi \bigl((b_0-b_0^m)\Delta_2 \bigr) Q \Bigr)
\end{equation*}
for $m \ge 2$. This is also affirmative at $m=1$. Hence,
\begin{equation*}
\exp (2\pi i H)-1 = \mathrm{diag} \Bigl( \bigl( \exp (2\pi i \phi (a_0\Delta_1 ))-1 \bigr) Q,\ \bigl( \exp (2\pi i \phi (b_0\Delta_2 ))-1 \bigr) Q \Bigr).
\end{equation*}
Consequently,
\begin{eqnarray*}
\exp (2\pi i H)=
\begin{pmatrix}
\exp (2\pi i \phi (a_0\Delta_1 ))Q & 0\\
0 & \exp (2\pi i \phi (b_0\Delta_2 ))Q \\
\end{pmatrix}
+\begin{pmatrix}
1-Q & 0 \\
0 & 1-Q \\
\end{pmatrix}.
\end{eqnarray*}
Next, the map $K_1(A) \longrightarrow K_1(K(F(X_{(1,2)})))$ defined by
\begin{equation*}
K_1(A) \ni [z]_1 \longmapsto [zQ+(1-Q)]_1 \in K_1(K(F(X_{(1,2)})))
\end{equation*}
is a group isomorphism (cf. \cite{Pimsner}): hence, $\delta_0[P]_0=[\exp (2\pi i H)]_1$ can be regarded as 
\begin{equation*}
\bigl[\ \mathrm{diag}\bigl( \exp (2\pi i \phi (a_0\Delta_1 )) ,\ \exp (2\pi i \phi (b_0\Delta_2 )) \bigr)\ \bigr]_1 \mbox{\quad in } K_1(C(\mathbb{T} )).
\end{equation*}
From the graphs of $a_0,b_0$, this element is $[z^{-4}]_1$ in $K_1(C(\mathbb{T} ))$. Hence, we have finished the proof of $[P]_0=-4$.
\end{proof}

\section{Subalgebras of $\mathcal{O}_{(m,n)}(\mathbb{T}) $ and cyclic group actions}
\label{sec:Subalgebra}

\subsection{Subalgebras of $\mathcal{O}_{(m,n)}(\mathbb{T}) $}

Let $\theta$ be an irrational number. Then the $C^*$-subalgebra $C^*(u^k,v)$ of the irrational rotation algebra $A_\theta$ generated by $u^k$ and $v$ is isomorphic to $A_{k \theta}$ from the relation $u^kv=e^{2\pi i k \theta}vu^k$ and from simplicity. Moreover, the subalgebra $C^*(u,v^k)$ of $A_\theta$ is isomorphic to $A_{k\theta}$. We would like to consider the corresponding problem for $\mathcal{O}_{(m,n)}(\mathbb{T}) $ i.e., we shall discuss the $C^*$-subalgebras $C^*(z^k,S_1)$ and $C^*(z,S_1^k)$ of $\mathcal{O}_{(m,n)}(\mathbb{T}) $. First, we consider some easy cases.

\begin{lem}\normalfont\slshape
Consider $\mathcal{O}_{(m,n)}(\mathbb{T})$ for $m \ge 1, n \ge 2$ and $\mathrm{gcd}(m,n)=1$. Then we have the following:
\begin{enumerate}
\item{The $C^*$-subalgebra $C^*(z,S_1^k)$ of $\mathcal{O}_{(m,n)}(\mathbb{T})$ generated by $z$ and $S_1^k$ is isomorphic to $\mathcal{O}_{(m^k,n^k)}(\mathbb{T})$.}
\item{If $k|n$ ($n$ is divided by $k$), then the $C^*$-subalgebra $C^*(z^k,S_1)$ of $\mathcal{O}_{(m,n)}(\mathbb{T})$ generated by $z^k$ and $S_1$ is isomorphic to $\mathcal{O}_{(m,n)}(\mathbb{T}).$}
\end{enumerate}
\end{lem}

\begin{proof}
(1) Put $\widetilde{S}_j=z^{j-1}S_1^k$ for $j=1,\cdots ,n^k$. Then we can easily check the relations
\begin{equation*}
z\widetilde{S}_{n^k}=\widetilde{S}_1z^{m^k}, \quad \widetilde{S}_i^*\widetilde{S}_j=\delta_{ij} ,\quad \sum_{i=1}^{n^k} \widetilde{S}_i\widetilde{S}_i^*=1.
\end{equation*}
Hence, $\mathcal{O}_{(m,n)}(\mathbb{T}) \cong C^*(z,S_1^k)$.\\
(2) Since $k|n$, there is a $p \in \mathbb{N}$ such that $n=kl$. Hence, $z^n \in C^*(z^k,S_1)$. The relation $S_1^*z^nS_1=z^m$ implies that $z^m \in C^*(z^k,S_1)$. Since $\mathrm{gcd}(m,n)=1$, there exist $p,q \in \mathbb{Z}$ such that $mp+nq=1$ and hence $z \in C^*(z^k,S_1)$. Since $z$ is one of the generators of $\mathcal{O}_{(m,n)}(\mathbb{T}) $, we have $C^*(z^k,S_1) \cong \mathcal{O}_{(m,n)}(\mathbb{T}) $.
\end{proof}

Below, we discuss for the general case of $k \in \mathbb{N}$.

\begin{prop}\normalfont\slshape
\label{prop:isom}
For $m \ge 1, n \ge 2$, and $\mathrm{gcd}(m,n)=1$, consider a $C^*$-algebra $\mathcal{O}_{(m,n)}$ generated by $z$ and $S_1$. Then the $C^*$-subalgebra $C^*(z^k,S_1)$ of $\mathcal{O}_{(m,n)}(\mathbb{T})$ generated by $z^k$ and $S_1$ is isomorphic to $\mathcal{O}_{(m,n)}(\mathbb{T}) $ for any $k \in \mathbb{N}$.
\end{prop}

\begin{proof}
 If $\mathrm{gcd}(k,n)\ge 2$, then we can reduce the case of $\mathrm{gcd}(k,n)=1$ using the relation $S_1^*z^nS_1=z^m$. Let represent $k,n$ as $k=p^\alpha p_1, n=p^\beta p_2$, where $p$ is a prime number and $p_1,p_2$ do not contain $p$ as factors. Suppose that $\alpha >\beta$ then $kp_2=(p^{\alpha - \beta }p_1)n$ and
\begin{equation*}
C^*(z^k,S_1) \ni S_1^*z^{kp_2}S_1=S_1^* z^{(p^{\alpha - \beta }p_1)n}S_1=z^{p^{\alpha - \beta }p_1m}.
\end{equation*}
Since $z^{p^{\alpha - \beta }p_1m}, z^{p^{\alpha - \beta }p_1n} \in C^*(z^k,S_1)$ and $\mathrm{gcd}(m,n)=1$, we obtain $z^{(p^{\alpha - \beta }p_1)} \in C^*(z^k,S_1)$. Repeating this process, we may assume that $\alpha \le \beta$. Then $kp^{\beta - \alpha}p_2=p_1n$ and
\begin{equation*}
C^*(z^k,S_1) \ni S_1^*z^{kp^{\beta - \alpha}p_2}S_1=S_1^*z^{p_1n}S_1=z^{p_1m}.
\end{equation*}
which also implies that $z^{p_1} \in C^*(z^k,S_1)$. Moreover $C^*(z^k,S_1)=C^*(z^{p_1},S_1)$ because $(z^{p_1})^{p^\alpha}=z^k$. Continuing this process, we can reduce the case of $\mathrm{gcd}(k,n)=1$.
Hence, it is enough to consider the case of $\mathrm{gcd}(k,n)=1$.

 Since $\mathrm{gcd}(k,n)=1$,
\begin{equation*}
\{ (q-1)k\ (\mbox{mod}\ n) | 1 \le q \le n\} =\mathbb{Z}_n.
\end{equation*}
Hence, for any $1 \le q \le n$, there exists $0 \le l_q \le (n-1)$ and $p_q \in \mathbb{Z}$ such that $(q-1)k=l_q + np_q$. For $1 \le q \le n$, put $\widetilde{S}_q=z^{(q-1)k}S_1$. Then
\begin{eqnarray*}
\sum_{q=1}^{n} \widetilde{S}_q \widetilde{S}_q^*
&=&\sum_{q=1}^{n} (z^{qk}S_1)(z^{qk}S_1)^*
=\sum_{q=1}^{n} (z^{l_q}S_1z^{mp_q})(z^{l_q}S_1z^{mp_q})^*
=\sum_{q=1}^{n} (z^{l_q}S_1)(z^{l_q}S_1)^*\\
&=&\sum_{i=1}^{n} S_iS_i^* =1
\end{eqnarray*}
If we put $w:=z^k$, then $w$ is a full spectrum unitary and for $1 \le i \le n-1$,
\begin{equation*}
w\widetilde{S}_i =\widetilde{S}_{i+1} , \quad w\widetilde{S}_{n}=\widetilde{S}_1w^m , \quad \sum_{k=1}^{n}\widetilde{S}_k \widetilde{S}_k^*=1.
\end{equation*}
Hence, $\mathcal{O}_{(m,n)}(\mathbb{T})  \cong C^*(z^k,S_1)=C^*(w,\widetilde{S}_1)$
\end{proof}


\subsection{Actions of cyclic groups on $\mathcal{O}_{(m,n)}(\mathbb{T})$}
\label{sec:action}

In \cite{IKW}, for the $C^*$-algebras associated complex dynamical systems, Izumi-Kajiwara-Watatani studied automorphisms arising from symmetries of the dynamical systems. In the case of $R(z)=z^n$ (whose $C^*$-algebra is $\mathcal{O}_{(1,n)}$), the dihedral group $\mathbb{Z}_{n-1} \rtimes \mathbb{Z}_{2}$ acts on the $C^*$-algebra $\mathcal{O}_{(1,n)}(\mathbb{T}) $ (see Example 8.3 in \cite{IKW}). They show that this dihedral action is outer. In this subsection, we shall consider extending this action on $\mathcal{O}_{(m,n)}(\mathbb{T} )$. Define a $\mathbb{Z}_{|n-m|}$-action $\beta : \mathbb{Z}_{|n-m|} \longrightarrow \mathrm{Aut}(\mathcal{O}_{(m,n)}(\mathbb{T}))$ by
\begin{equation*}
\beta_t : \quad z \longmapsto tz ,\quad S_1 \longmapsto S_1.
\end{equation*}
for $t \in \mathbb{Z}_{|n-m|}$. Furthermore, we shall define a $\mathbb{Z}_2$-action $\sigma\in \mathrm{Aut}(\mathcal{O}_{(m,n)}(\mathbb{T}))$ by
\begin{equation*}
\sigma : \quad z \longmapsto z^{-1} ,\quad S_1 \longmapsto S_1.
\end{equation*}
We can easily check that these actions are well-defined. We shall show that these actions are outer by an elementary method.

\begin{prop}\normalfont\slshape
\label{prop:outer}
Suppose that $m \ge 1, n \ge 2$ and $\mathrm{gcd}(m,n)=1$. Let $\beta$, $\sigma$ be group actions on $\mathcal{O}_{(m,n)}(\mathbb{T})$ defined as above.\\
If $|n-m| \ge 2$, then $\mathbb{Z}_{|n-m|}$-action $\beta$ is outer. The $\mathbb{Z}_2$-action $\sigma$ is also outer.
\end{prop}

\begin{proof}
First, we shall show that $\beta$ is outer for $|n-m| \ge 2$. This is an well-known argument, but we give a proof for completeness. Suppose that there exists $t \in \mathbb{Z}_{|n-m|}$ with $t \ne 1$ and $u \in \mathcal{O}_{(m,n)}(\mathbb{T})$ such that $\beta_t(x)=uxu^*$. Let us represent $\mathcal{O}_{(m,n)}(\mathbb{T})$ on the Hilbert space $l^2(\mathbb{Z}[1/m])$ with CONS $\{ e_k\}_{k \in \mathbb{Z}[1/m]}$ by
\begin{equation*}
ze_k=e_{k+1} ,\quad S_1e_k=e_{\frac{n}{m}k-1}.
\end{equation*}
We can check that this representation is well-defined. In particular, $S_2e_0=e_0$. Since $\beta_t(S_2)=tS_2$, we get $tu^*e_0=u^*e_0$, which implies that $t=1$. This is a contradiction, so we conclude that $\beta$ is outer for $|n-m| \ge 2$.

Next, we shall show that $\sigma$ is outer. Suppose that there is a unitary $w \in \mathcal{O}_{(m,n)}(\mathbb{T})$ such that $\sigma (x)=wxw^*$ for $x \in \mathcal{O}_{(m,n)}(\mathbb{T})$. We consider another representation on $l^2(\mathbb{Z}[1/m])$ by
\begin{equation*}
ze_k=e_{k+1} ,\quad S_1e_k=e_{\frac{n}{m}k}.
\end{equation*}
Then we shall show that $w$ can be identified by $we_k=e_{-k}$ ($k \in \mathbb{Z}[1/m]$). Since $wS_1=S_1w$, 
\begin{equation*}
we_0=wS_1^ke_0=S_1^kwe_0 \in S_1^kl^2(\mathbb{Z}[1/m] ) = \overline{\mathrm{span}}\{ \cdots ,e_{-(n/m)^k},e_0,e_{(n/m)^k} \cdots \} 
\end{equation*}
for any $k \in \mathbb{N}$. This says that there exists $\lambda \in \mathbb{T}$ such that $we_0=\lambda e_0$ since $w$ is a unitary. Because $wz^k= z^{-k}w$ for any $k \in \mathbb{Z}$,
\begin{equation*}
we_k=wz^ke_0=z^{-k}we_0=z^{-k}\lambda e_0= \lambda e_{-k} .
\end{equation*}
For $p/m^q \in \mathbb{Z}[1/m]$ $(p \in \mathbb{Z} ,q \in \mathbb{N}_0)$, there exists $r,s \in \mathbb{Z}$ such that $p/m^q=(n/m)^qr+s$. Then
\begin{equation*}
we_{p/m^q}=wz^s S_1^q e_r=z^{-s} S_1^q w e_r=\lambda z^{-s} S_1^q e_{-r}=\lambda e_{-(n/m)^qr-s}=\lambda e_{-p/m^q}.
\end{equation*}
Hence, we can assume that $we_k=e_{-k}$ for $k \in \mathbb{Z}[1/m]$.\\
We would like to obtain a contradiction for $w \in \mathcal{O}_{(m,n)}(\mathbb{T})$. Since we supposed that $w \in \mathcal{O}_{(m,n)}(\mathbb{T})$, there exists $K \in \mathbb{N}$ , $\mu_i,\nu_i \in \mathcal{W}, k_i \in \mathbb{Z} $ such that
\begin{equation*}
\Bigl\| w- \sum_{i=1}^K S_{\mu_i}z^{k_i}S_{\nu_i}^* \Bigr\| < 1.
\end{equation*}
On the other hand, for each monomial $S_{\mu_i}z^{k_i}S_{\nu_i}^*$, there exists at most one element $p_i \in \mathbb{Z}[1/m]$ such that $S_{\mu_i}z^{k_i}S_{\nu_i}^*e_{p_i}=e_{-p_i}$. Hence for $q \in \mathbb{Z}[1/m] \setminus \{p_i \}_{i=1}^K$,
\begin{equation*}
\Bigl\| w- \sum_{i=1}^K S_{\mu_i}z^{k_i}S_{\nu_i}^* \Bigr\| 
\ge \Bigl\| \theta_{-q, -q} \Bigl( w- \sum_{i=1}^K S_{\mu_i}z^{k_i}S_{\nu_i}^* \Bigr) e_q \Bigr\| 
\ge \| e_{-q} \| =1
\end{equation*}
where $\theta_{i,j}(\xi )= \langle \xi, e_i \rangle_{l^2(\mathbb{Z}[1/m])}e_j$ is a rank-one operator. This is a contradiction. Hence $\sigma$ is outer.
\end{proof}


\subsection{Fixed point algebra of $\mathbb{Z}_{|n-m|}$-action}

We shall check the fixed point algebra of $\mathbb{Z}_{|n-m|}$-action $\beta$. Let us denote $k \equiv l$ by $k=l \ \mathrm{mod}\ |n-m|$.

\begin{prop}\normalfont\slshape
Suppose that $m \ge 1, n \ge 2$, and $\mathrm{gcd}(m,n)=1$. Moreover, we assume that $|n-m| \ge 2$. Then the fixed point algebra of $\mathbb{Z}_{|n-m|}$-action $\beta$ defined in Section \ref{sec:action} is the $C^*$-subalgebra $C^*(z^{|n-m|},S_1)$ of $\mathcal{O}_{(m,n)}(\mathbb{T})$ generated by $z^{|n-m|}$ and $S_1$, which is isomorphic to $\mathcal{O}_{(m,n)}(\mathbb{T}) $ by Proposition \ref{prop:isom}.
\end{prop}

\begin{proof}
The inclusion $C^*(z^{|n-m|},S_1) \subset \mathcal{O}_{(m,n)}(\mathbb{T})^\beta$ is trivial. We shall show that $\mathcal{O}_{(m,n)}(\mathbb{T})^\beta \subset C^*(z^{|n-m|},S_1)$.\\
Since $\mathcal{O}_{(m,n)}(\mathbb{T})$ is spanned by $S_\mu z^k S_\nu^*$ $(k \in \mathbb{Z}, \mu ,\nu \in \mathcal{W}_n)$, we shall consider a monomial $S_\mu z^k S_\nu^*$. Suppose that $S_\mu z^k S_\nu^* \in \mathcal{O}_{(m,n)}(\mathbb{T})^\beta$: then the indices $k,\mu ,\nu$ satisfy
\begin{equation*}
\sum_{i=1}^{|\mu |} (\mu_i-1) + k - \sum_{j=1}^{|\nu |} (\nu_j-1) \equiv 0.
\end{equation*}
Since $\mathrm{gcd}(n,|n-m|)=1$, there is an integer $p_1$ such that $\gamma_1:=(\mu_1-1)+p_1n \equiv 0$. Moreover, for $\mu_2-1-p_1m \in \mathbb{Z}$, there is an integer $p_2$ such that $\gamma_2:=(\mu_2-1-p_1m)+p_2n \equiv 0$. Repeating this process, we find that there are $\{ p_i \}_{i=2}^{| \mu |}$ such that $\gamma_i := (\mu_i-1-p_{i-1}m)+p_in \equiv 0$ $(i=2,\cdots ,|\mu |)$. Hence,
\begin{eqnarray*}
S_\mu
&=&z^{\mu_1-1}S_1 z^{\mu_2-1} S_1 \cdots z^{\mu_{|\mu |}-1} S_1\\
&=&z^{\mu_1-1+p_1n}S_1 z^{(\mu_2-1-p_1m)+p_2n}S_1 \cdots z^{(\mu_{|\mu|}-1-p_{|\mu|-1}m + p_{|\mu |n}}S_1 z^{-p_{|\mu |}m}\\
&=&z^{\gamma_1}S_1 z^{\gamma_2}S_1 \cdots z^{\gamma_{|\mu |}}S_1 z^{-p_{|\mu |}m}
\end{eqnarray*}
Furthermore, there are integers $\{ q_i \}_{i=1}^{|\nu |}$ such that $\delta_1:=(\nu_1-1)+p_1n \equiv 0$ and $\delta_i := (\nu_i-1-p_{i-1}m)+p_in \equiv 0$ $(i=2,\cdots ,|\nu |)$. Then
\begin{equation*}
S_\nu=z^{\delta_1}S_1 z^{\delta_2}S_1 \cdots z^{\delta_{|\nu |}}S_1 z^{-q_{|\nu |}m}.
\end{equation*}
Hence, we obtain
\begin{equation*}
S_\mu z^k S_\nu^* = z^{\gamma_1}S_1 z^{\gamma_2}S_1 \cdots z^{\gamma_{|\mu |}}S_1 z^{-p_{|\mu |}m+k+q_{| \nu |}m}
S_1^* z^{-\delta_{|\nu |}} \cdots S_1^* z^{-\delta_2} S_1^*z^{-\delta_1}.
\end{equation*}
Then
\begin{eqnarray*}
0 &\equiv& 
\sum_{i=1}^{|\mu |} \gamma_i - \sum_{i=1}^{|\nu |} \delta_i \\
&=&\Bigl( \sum_{i=1}^{|\mu |} (\mu_i-1)-\sum_{i=1}^{|\nu |} (\nu_i-1) \Bigr) + \sum_{i=1}^{|\mu |-1}p_i(n-m) - \sum_{i=1}^{|\nu |-1} q_i(n-m) + p_{| \mu |}n - q_{| \nu |}n\\
&\equiv&-k+p_{| \mu |}m - q_{| \nu |}m.
\end{eqnarray*}
Hence, we obtain $\gamma_i , \delta_j , -p_{|\mu |}m+k+q_{| \nu |}m \in |n-m| \mathbb{Z}$ and this implies that $S_\mu z^k S_\nu \in C^*(z^{|n-m|},S_1)$. Hence we conclude that $\mathcal{O}_{(m,n)}(\mathbb{T})^\beta \subset C^*(z^{|n-m|},S_1)$.
\end{proof}

\subsection{Fixed point algebra of symmetric action}
\label{sec:fixsym}

We shall determine the fixed point algebra of the $\mathbb{Z}_2$-action $\sigma$ defined in Section \ref{sec:action}. In this subsection, we shall show the following proposition:

\begin{prop}\normalfont\slshape
\label{prop:symmetry}
Suppose that $m \ge 1, n \ge 2 $, and $\mathrm{gcd}(m,n)=1$. Let $\sigma$ be the $\mathbb{Z}_2$-action on $\mathcal{O}_{(m,n)}(\mathbb{T})$ defined in Section \ref{sec:action}. Then the fixed point algebra $\mathcal{O}_{(m,n)}(\mathbb{T})^\sigma $ is a Kirchberg algebra satisfying UCT and its $K$-group is following: if $m$ is even, then
\begin{equation*}
K_0(\mathcal{O}_{(m,n)}(\mathbb{T})^\sigma )=\mathbb{Z}_{n-1}, \quad K_1(\mathcal{O}_{(m,n)}(\mathbb{T})^\sigma )=0 ,\quad [1]_0=0
\end{equation*}
and if $m$ is odd, then
\begin{equation*}
K_0(\mathcal{O}_{(m,n)}(\mathbb{T})^\sigma )=\mathbb{Z}[1]_0 \oplus \mathbb{Z}_{n-1} \oplus \mathbb{Z}_{n-1}, \quad K_1(\mathcal{O}_{(m,n)}(\mathbb{T})^\sigma )=\mathbb{Z}.
\end{equation*}
\end{prop}

\begin{proof}
First, we compute the $K$-groups of $\mathcal{O}_{(m,n)}(\mathbb{T}) ^\sigma$. Since $\mathcal{O}_{(m,n)}(\mathbb{T})$ is simple and $\sigma$ is outer, $\mathcal{O}_{(m,n)}(\mathbb{T}) ^\sigma$ is Morita equivalent to $\mathcal{O}_{(m,n)}(\mathbb{T})  \rtimes_\sigma \mathbb{Z}_2$. Hence the group $K_i(\mathcal{O}_{(m,n)}(\mathbb{T}) ^\sigma )$ is isomorphic to $K_i(\mathcal{O}_{(m,n)}(\mathbb{T})  \rtimes_\sigma \mathbb{Z}_2)$ ($i=0,1$). If $w$ is the unitary of $\mathcal{O}_{(m,n)}(\mathbb{T})  \rtimes_\sigma \mathbb{Z}_2$ implementing $\sigma$, then elements in the crossed product have the form $x+yw$ where $x,y \in \mathcal{O}_{(m,n)}(\mathbb{T}) $. To compute $K$-groups, we construct another Cuntz-Pimsner algebra. We define a $\mathbb{Z}_2$-action $\sigma_0$ on $A_0=C(\mathbb{T} )$ by $\sigma_0(z_0)=z_0^{-1}$ (where $z_0$ is the unitary generator of $A_0$) and define $B=A_0 \rtimes_{\sigma_0} \mathbb{Z}_2$ with implement unitary $w_0$. Let $Y=B^{\oplus n}$ with the right $B$-action defined by $(b_i)_{i=1}^n \cdot b=(b_ib)_{i=1}^n$ for $(b_i)_{i=1}^n \in Y, b \in B$. Let us put $u_i'=(0,\cdots ,1,\cdots ,0) \in Y$ for $i=1 ,\cdots , n$. Define a $B$-valued inner product by
\begin{equation*}
\langle (b_i)_i , (b_i')_i \rangle_B=\sum_{i=1}^n b_i^*b_i'.
\end{equation*}
Then $Y$ is a full Hilbert $B$-module with this inner product. Let us define a left $B$-action $\phi:B \longrightarrow L_B(Y)\cong M_n(B)$ by
\[
\phi(z_0)=
\left(
\begin{array}{cccc|c}
0&0&\cdots&0&z_0^m\\
\hline
1&0&\cdots&0&0\\
0&1&\cdots&0&0\\
\vdots&\vdots&\ddots&\vdots&\vdots\\
0&0&\cdots&1&0
\end{array}
\right)
,\ 
\phi(w_0)=
\left(
\begin{array}{c|cccc}
w_0&0&\cdots&0&0\\
\hline
0&0&\cdots&0&w_0z_0^m\\
0&0&\cdots&w_0z_0^m&0\\
\vdots&\vdots&\ddots&\vdots&\vdots\\
0&w_0z_0^m&\cdots&0&0
\end{array}
\right)
\]
Then we can check that $\phi (w_0)$ is self-adjoint unitary and that $\phi(w_0) \phi(z_0) \phi(w_0) = \phi (z_0)^{-1}$. Hence, $\phi$ is a *-homomorphism. In fact, $\phi$ is faithful:

\begin{lem}\normalfont\slshape
$\phi$ is faithful.
\end{lem}

\begin{proof}
Define $E$ to be the canonical faithful conditional expectation from $B$ onto $A_0$. Let $\widehat{\sigma_0}$ be the dual action of $\sigma_0$ and let $u$ be the implement unitary of $B \rtimes_{\widehat{\sigma_0}} \mathbb{Z}_2$. Define a unitary in $M_n(B \rtimes_{\widehat{\sigma_0}} \mathbb{Z}_2)$ by $U=\mathrm{diag}(u,\cdots , u )$. Then we can check that $\mathrm{Ad}(U)(\phi (z_0)):=U\phi(z_0)U^*=\phi(z_0)$ and $\mathrm{Ad}(U)(\phi (w_0))=-\phi (w_0)$. Set $E_1:=\frac{1}{2}(\mathrm{id}+\mathrm{Ad}(U) ) : \phi(B) \longrightarrow \phi (B)$, Then $E_1$ is a faithful conditional expectation onto $\phi (A_0)$ and satisfies $E_1 \circ \phi = \phi \circ E$. Since $\phi |_{A_0} :A_0 \longrightarrow \phi (A_0)$ is an isomorphism, the equation $E_1 \circ \phi = \phi \circ E$ induces the faithfulness of $\phi$.
\end{proof}
Hence, we can construct the Cuntz-Pimsner algebra $\mathcal{O}_Y$ from these data. Then $\mathcal{O}_{(m,n)}(\mathbb{T})  \rtimes_\sigma \mathbb{Z}_2$ is isomorphic to $\mathcal{O}_Y$ by both universalities; the isomorphism is determined by
\begin{equation*}
\psi (z)=z_0 ,\ \psi(w)=w_0 ,\ \psi(S_i)=S_{u_i'} \ (i=1,\cdots n).
\end{equation*}
Next, let us compute the $K$-groups of $\mathcal{O}_Y$ instead of those of $\mathcal{O}_{(m,n)}(\mathbb{T})  \rtimes_\sigma \mathbb{Z}_2$. From the six-term exact sequence of Cuntz-Pimsner algebra, we obtain the exact sequence
\[
\begin{CD}
K_0(B) @>\mathrm{id}_*-[Y]_0>> K_0(B) @>\iota_*>> K_0(\mathcal{O}_Y )\\
@A\delta_1AA @. @VV\delta_0V \\
K_1(\mathcal{O}_Y) @<\iota_*<< K_1(B) @<\mathrm{id}_*-[Y]_1<< K_1(B)
\end{CD}
\]
where $[Y]_i\ (i=0,1)$ is the group homomorphism arising from $\phi$.

We recall that $K_0(B)=\mathbb{Z}^3$ with generators $e_0:=[1]_0, e_1:=[\frac{1}{2}(1+w_0)]_0, e_2:=[\frac{1}{2}(1+w_0z_0)]_0$ and $K_1(B)=0$. We can show
\begin{eqnarray*}
[Y]_0(e_0)=ne_0 ,
\quad [Y]_0(e_1)=e_1+(n-1)\Bigl[\frac{1}{2}(1+w_0z_0^m)\Bigr]_0 ,\quad 
[Y]_0(e_2)=n\Bigl[\frac{1}{2}(1+w_0z_0^m)\Bigr]_0
\end{eqnarray*}
(see the proof of Theorem \ref{thm:Katsura}). 
Since
\[
\Bigl[\frac{1}{2}(1+w_0z_0^m)\Bigr]_0=
\left\{ 
 \begin{array}{ll}
   \bigl[\frac{1}{2}(1+w_0)\bigl]_0 & m \mathrm{:even} \\
   \bigl[\frac{1}{2}(1+w_0z_0)\bigr]_0 & m \mathrm{:odd}, \\
 \end{array}
\right.
\]
if $m$ is even, then $[Y]_0(e_1)=ne_1$, $[Y]_0(e_2)=ne_1$ and if $m$ is odd, then $[Y]_0(e_1)=e_1+(n-1)e_2$ and $[Y]_0(e_2)=ne_2$. From the six-term exact sequence, we can conclude that
\[
K_0(\mathcal{O}_Y)=
\left\{ 
 \begin{array}{ll}
   \mathbb{Z}_{n-1}e_0 & m \mathrm{:even} \\
   \mathbb{Z}_{n-1}e_0 \oplus \mathbb{Z}e_1 \oplus \mathbb{Z}_{n-1}e_2 & m \mathrm{:odd} \\
 \end{array}
\right.
,\quad
K_1(\mathcal{O}_Y)=
\left\{ 
 \begin{array}{ll}
   0 & m \mathrm{:even} \\
   \mathbb{Z} & m \mathrm{:odd}, \\
 \end{array}
\right.
\]
Since the isomorphism from $K_0(\mathcal{O}_{(m,n)}(\mathbb{T})^\sigma )$ to $K_0(\mathcal{O}_{(m,n)}(\mathbb{T}) \rtimes_\sigma \mathbb{Z}_2)$ is defined by $[p]_0 \longmapsto [\frac{1}{2}(p+pw)]_0$, we have determined the K-groups of $\mathcal{O}_{(m,n)}(\mathbb{T})^\sigma$ in Proposition \ref{prop:symmetry}.

Since $\sigma$ is outer (Proposition \ref{prop:outer}) and $\mathcal{O}_{(m,n)}(\mathbb{T}) $ is purely infinite simple, $\mathcal{O}_{(m,n)}(\mathbb{T}) ^\sigma $ is purely infinite simple (by using Lemma 10 of \cite{KK} and $\mathcal{O}_{(m,n)}(\mathbb{T}) ^\sigma$ is a hereditary algebra of $\mathcal{O}_{(m,n)}(\mathbb{T})  \rtimes_\sigma \mathbb{Z}_2$). Since $\mathcal{O}_{(m,n)}(\mathbb{T}) ^\sigma$ is a hereditary algebra of $\mathcal{O}_{(m,n)}(\mathbb{T})  \rtimes_\sigma \mathbb{Z}_2$ and $\mathcal{O}_{(m,n)}(\mathbb{T})  \rtimes_\sigma \mathbb{Z}_2$ is nuclear, $\mathcal{O}_{(m,n)}(\mathbb{T}) ^\sigma$ is also nuclear. The separability of $\mathcal{O}_{(m,n)}(\mathbb{T}) ^\sigma$ is trivial. Moreover $\mathcal{O}_Y$ satisfies UCT because $B$ satisfies UCT, and $\mathcal{O}_{(m,n)}(\mathbb{T}) ^\sigma $ is Morita equivalent to $\mathcal{O}_Y$, so $\mathcal{O}_{(m,n)}(\mathbb{T}) ^\sigma $ also satisfies UCT. Hence, we have completed the proof of Proposition \ref{prop:symmetry}.

\end{proof}


\section{Entropy estimate for the canonical endomorphism on $\mathcal{O}_{(1,n)}(\mathbb{T}) $}
\label{sec:Entropy}

For the Cuntz algebra $\mathcal{O}_n$, Choda has computed Voiculescu's entropy (\cite{Choda}) for the \textit{Cuntz's canonical endomorphism} defined by
\begin{equation*}
\Phi_0 (x)=\sum_{i=1}^n S_ixS_i^* , \quad (x \in \mathcal{O}_n ).
\end{equation*}
In this section, we consider an analogy to this problem for $\mathcal{O}_{(1,n)}(\mathbb{T})$

For $\mathrm{gcd}(m,n)=1$ and $n \ge 2$, let us defined the \textit{canonical endomorphism} on $\mathcal{O}_{(m,n)}$ by
\begin{equation*}
\Phi (x)=\sum_{i=1}^n S_ixS_i^* , \quad (x \in \mathcal{O}_{(m,n)}(\mathbb{T}) ).
\end{equation*}
Its name is derived from one of the Cuntz algebras $\mathcal{O}_n$. Then $\Phi$ is a *-endomorphism on $\mathcal{O}_{(m,n)}$.
We would like to compute Voiculescu's topological entropy for $\Phi$ on $\mathcal{O}_{(1,n)}(\mathbb{T})$. Our method is similar to that of Boca-Goldstein(\cite{BG})\\
Let us recall the definition of the Voiculescu's topological entropy(\cite{Voi}). Let $B$ be a nuclear $C^*$-algebra with unit. Let $CPA(B)$ be the triples $(\phi , \psi , B)$, where $C$ is a finite-dimensional $C^*$-algebra, and $\phi : B \longrightarrow C$ and $\psi :C \longrightarrow B$ are unital completely positive maps. Let $Pf(B)$ be the set of finite subsets of $B$. For an $\omega \in Pf(B)$, put
\begin{equation*}
rcp(\omega ;\delta )=\inf \{ \mathrm{rank}\ C : (\phi , \psi , C) \in CPA(B) , \| \psi \circ \phi (a)-a \| <\delta ,a \in B \}
\end{equation*}
where $\mathrm{rank}\ C $ means the dimension of a maximal abelian self-adjoint subalgebra of $C$. Since $B$ is nuclear, for any $\omega \in Pf(B)$ and $\delta >0$, there exists $(\phi , \psi , C) \in CPA(B)$ such that $\| \psi \circ \phi (a)-a \| <\delta , a \in \omega $. For a unital *-endomorphism $\beta$ of $B$, put
\begin{equation*}
ht(\beta ,\omega ;\delta )=\limsup_{N \longrightarrow \infty}
\frac{1}{N} \log rcp (\omega \cup \beta (\omega ) \cup \cdots \cup \beta^{N-1}(\omega ) ;\delta )
\end{equation*}
and 
\begin{equation*}
ht(\beta ;\omega )=\sup_{\delta >0} ht(\beta , \omega ;\delta ).
\end{equation*}
Then \textit{(Voiculescu's) topological entropy} $ht(\beta )$ of $\beta$ is defined by
\begin{equation*}
ht(\beta)=\sup_{\omega \in Pf(B)} ht(\beta , \omega ). 
\end{equation*}
We recall the Kolmogorov-Sinal type theorem.

\begin{thm}[Voiculescu \cite{Voi}]\normalfont\slshape
\label{thm:KSthm}
Let $\omega_j \in Pf(B)$ such that $\omega_1 \subset \omega_2 \subset \cdots $ and the linear span of $\bigcup_{j \in \mathbb{N} } \omega_j$ is dense in $B$. Then
\begin{equation*}
ht(\beta) = \sup_{j \in \mathbb{N}} ht(\beta , \omega_j).
\end{equation*}
\end{thm}

Let $\varphi$ be a state of $B$ with $\varphi \circ \beta =\varphi$. An estimate between $ht(\beta )$ and \textit{Connes-Narnhofer-Thirring (CNT) entropy} $h_\varphi (\beta)$ (\cite{CNT}) is given by
\begin{equation*}
h_\varphi (\beta ) \leqq ht(\beta ).
\end{equation*}
which was proved by Voiculescu(\cite{Voi}).

The $C^*$-algebra $\mathcal{O}_{(1,n)}(\mathbb{T})$ has exactly one $\log n$-KMS state $\varphi$ for the gauge action of $\mathcal{O}_{(1,n)}$. This KMS-state is written as $\varphi=\tau \circ E$, where $\tau$ is the unique normalized trace on the $(1,n)$-type Bunce-Deddens algebra and $E$ is the conditional expectation onto the Bunce-Deddens algebra. Our main theorem in this section is as follows:

\begin{thm}\normalfont\slshape
Suppose $m=1$ and $n \ge 2$. Let $\varphi$ be the unique $\log n$-KMS state for gauge action of $\mathcal{O}_{(1,n)}(\mathbb{T}) $. Let $\Phi$ be the canonical endomorphism defined as above. Then the Voiculescu's topological entropy $ht(\Phi )$ for $\Phi$ and the CNT-entropy $h_\varphi (\Phi )$ for $\Phi$ and $\varphi$ are both equal to $\log n$;
\begin{equation*}
h_\varphi (\Phi )=ht(\Phi )=\log n.
\end{equation*}
\end{thm}

\begin{proof}
First, we define a map $\rho_r : \mathcal{O}_{(1,n)}(\mathbb{T})  \longrightarrow M_{n^r}(\mathbb{C}) \otimes \mathcal{O}_{(1,n)}(\mathbb{T}) $ for $r \ge 1$ by
\begin{equation*}
\rho_r(x)=\sum_{|\mu |,| \nu |=r} e_{\mu \nu} \otimes S_\mu^*xS_\nu .
\end{equation*}
We can check that this map is *-homomorphism and induce the isomorphism between $\mathcal{O}_{(1,n)}(\mathbb{T}) $ and $M_{n^r}(\mathbb{C}) \otimes \mathcal{O}_{(1,n)}(\mathbb{T}) $.\\
For $\mu \in \mathcal{W}_n^{(k)}$, we can see $\mu$ as $\sum_{i=1}^k (\mu_i-1)n^{i-1}$ via $S_\mu=S_{\mu_1}\cdots S_{\mu_k}$. Define 
\begin{equation*}
\| \mu \|_k=\sum_{i=1}^k (\mu_i-1) n^{i-1}
\end{equation*}
for $\mu \in \mathcal{W}_n^{(k)}$.
\begin{lem}\normalfont\slshape
Let $ N \ge 1$ and assume that $|\alpha |,| \beta | \le s$ and $N+s \le r$ and $1 \le l \le N$. Then for $|k| \le n^s$,
\[
\rho_r \circ \Phi^l(S_\alpha z^k S_\beta^*)
 =
\left\{
  \begin{array}{ll}
  x_0 \otimes z^{q_0}+x_1 \otimes z^{q_0+1} & |\alpha |=|\beta | \\
  \sum_{|\eta|=|\alpha |-|\beta |}
\Bigl( y_{0, \eta} \otimes S_{\eta}z^{q_0}+
y_{1, \eta } \otimes S_{\eta}z^{q_0+1} \Bigr) & | \alpha | > | \beta |,\ k \ge 0\\
  \sum_{|\eta|=|\alpha |-|\beta |}
\Bigl( y_{0, \eta} \otimes z^{q_0}S_{\eta}+
y_{1, \eta } \otimes z^{q_0+1}S_{\eta} \Bigr) & | \alpha | > | \beta |,\ k \le 0\\
  \sum_{|\eta|=|\alpha |-|\beta |}
\Bigl( y_{0, \eta} \otimes S_{\eta}^*z^{q_0}+
y_{1, \eta } \otimes S_{\eta}^*z^{q_0+1} \Bigr) & | \alpha | < | \beta |,\ k \ge 0\\
  \sum_{|\eta|=|\alpha |-|\beta |}
\Bigl( y_{0, \eta} \otimes z^{q_0}S_{\eta}^*+
y_{1, \eta } \otimes z^{q_0+1}S_{\eta}^* \Bigr) & | \alpha | < | \beta |,\ k \le 0\\

  \end{array}
\right.
\]
where $x_0,x_1,y_{0, \eta} ,y_{1, \eta }$ are partial isometries that depend on $\alpha ,\beta ,k$, and $|q_0| \le n^{s}$
\end{lem}

\begin{proof}
We consider the case of $k \ge 0$ (the case of $k \le 0$ is similar) and suppose that $|\beta | \le | \alpha |$.
\begin{eqnarray*}
\lefteqn{\rho_r\circ \Phi^l(S_\alpha z^kS_\beta^* )}\hspace{0cm}\\
&=&\sum_{|\mu | ,|\nu |=r} \sum_{| \gamma |=l} e_{\mu , \nu} \otimes S_\mu^* S_\gamma (S_\alpha z^k S_\beta^* ) S_\gamma^* S_\nu 
=\sum_{|\mu |=r-l-|\alpha | ,|\nu |=r-l-|\beta |} (\sum_{| \gamma |=l} e_{\mu \alpha \gamma , \nu \beta \gamma}) \otimes S_\mu^* z^k  S_\nu \\
\end{eqnarray*}
Let us put $x_{\mu , \nu}=\sum_{| \gamma |=l} e_{\mu \alpha \gamma , \nu \beta \gamma}$, and $p=r-l-|\beta |$. Then
\begin{eqnarray*}
\lefteqn{\rho_r\circ \Phi^l(S_\alpha z^kS_\beta^* )}\hspace{0cm}\\
&=&\sum_{|\mu |=r-l-|\alpha | ,|\nu |=p} x_{\mu , \nu} \otimes S_\mu^* z^k  S_\nu \\
&=&\sum_{|\mu |=r-l-|\alpha |}
\Bigl( \sum_{\{ \nu | 0 \le \| \nu \|_p \le n^p-k-1 \}} x_{\mu , \nu} \otimes S_\mu^* z^k  S_\nu
 + \sum_{q=1}^{n^{s-p}} \sum_{\{ \nu | qn^p-k \le \| \nu \|_p \le (q+1)n^p-k-1\}}x_{\mu , \nu} \otimes S_\mu^* z^k  S_\nu \Bigr)\\
&=&\sum_{|\mu |=r-l-|\alpha |}
\Bigl( \sum_{\{ \eta | k \le \| \eta \|_p \le n^p-1 \}} x_{\mu , \eta} \otimes S_\mu^* S_{\eta }
 + \sum_{q=1}^{n^{(s-p)}} \sum_{\{ \eta | qn^p \le \| \nu \|_p \le (q+1)n^p-1\}}x_{\mu , \eta} \otimes S_\mu^*  S_{\eta} z^q\Bigr)
\end{eqnarray*}
Let us define
\begin{equation*}
\mathcal{V}_0=\{ \eta \in \mathcal{W}_n^{(p)}| k \le \| \eta \|_p \le n^p-1 \},\quad 
\mathcal{V}_q=\{ \eta \in \mathcal{W}_n^{(p)}| qn^p \le \| \eta \|_p \le (q+1)n^p-1 \}
\end{equation*}
for $1 \le q \le n^{s-p}$. Then there exists $0 \le q_0 \le n^{(s-p)}$ such that all $\mathcal{V}_q$ are empty except $q=q_0$ or $q=q_0+1$. We shall define $x_{\mu,\eta}=0$ for $\| \eta \|_p < k, n-p+k \le \| \eta \|_p$. 
\begin{eqnarray*}
\rho_r\circ \Phi^l(S_\alpha z^kS_\beta^* )
=\sum_{|\mu |=r-l-|\alpha |}
\Bigl( \sum_{\eta \in \mathcal{V}_{q_0}} x_{\mu , \eta} \otimes S_\mu^* S_{\eta } z^{q_0}
 + \sum_{\eta \in \mathcal{V}_{q_0+1}}x_{\mu , \eta} \otimes S_\mu^*  S_{\eta} z^{q_0+1} \Bigr)
\end{eqnarray*}
If we assume $|\alpha |=|\beta |$, then
\begin{eqnarray*}
\rho_r\circ \Phi^l(S_\alpha z^kS_\beta^* )=
\Bigl(\sum_{|\mu |=r-l-|\alpha |}x_{\mu, \eta (\mu ) } \Bigr) \otimes z^{q_0}+
\Bigl(\sum_{|\mu |=r-l-|\alpha |}x_{\mu, \eta (\mu ) } \Bigr) \otimes z^{q_0+1},
\end{eqnarray*}
and $\sum_{|\mu |=r-l-|\alpha |}x_{\mu, \eta (\mu )}$ is a partial isometry. If $|\alpha |>|\beta |$, then
\begin{eqnarray*}
\rho_r\circ \Phi^l(S_\alpha z^kS_\beta^* )
=\sum_{|\eta|=|\alpha |-|\beta |}
\Bigl( x_{\mu (\eta), \eta} \otimes S_{\eta}z^{q_0}+
x_{\mu (\eta), \eta } \otimes S_{\eta}z^{q_0+1} \Bigr),
\end{eqnarray*}
where $x_{\mu (\eta), \eta}$ are partial isometries.\\
If we take the involution, we get the case of $| \alpha | < | \beta | $.
\end{proof}

Let us define
\begin{equation*}
\omega (s)=\{ S_\alpha z^k S_\beta^* |\ |\beta | , | \alpha | \le s, \ | k| \le n^{s} \}
\end{equation*}
which is increasing for $s \in \mathbb{N}$ and the linear span of the union is dense in $\mathcal{O}_{(1,n)}(\mathbb{T}) $. Since $\mathcal{O}_{(1,n)}(\mathbb{T}) $ is nuclear, there exist unital completely positive maps $\phi_0 : \mathcal{O}_{(1,n)}(\mathbb{T})  \longrightarrow M_{R}(\mathbb{C})$ and $\psi_0 : M_{R}(\mathbb{C}) \longrightarrow \mathcal{O}_{(1,n)}(\mathbb{T}) $ such that
\begin{eqnarray*}
\lefteqn{\sum_{|q| \le n^{s}} \Bigl( \| \psi_0 \circ \phi_0(S_\eta z^q )-S_\eta z^q\| + 
\| \psi_0 \circ \phi_0(z^qS_\eta  )-z^qS_\eta \| }\hspace{1cm}\\
& &+ \| \psi_0 \circ \phi_0(z^q S_\eta^* )- z^qS_\eta^*\| + 
\| \psi_0 \circ \phi_0(S_\eta^*z^q  )-S_\eta^*z^q\| \Bigr)
< \frac{\delta}{n^{s}}
\end{eqnarray*}
for $0 \le | \eta | \le s , |q| \le n^{s}$. Let us define $\phi : \mathcal{O}_{(1,n)}(\mathbb{T})  \longrightarrow M_{R}(\mathbb{C}) \otimes M_{n^r}(\mathbb{C})$ and $\psi : M_{R}(\mathbb{C}) \otimes M_{n^r}(\mathbb{C}) \longrightarrow \mathcal{O}_{(1,n)}(\mathbb{T}) $ by 
\begin{equation*}
\phi = (\mathrm{id} \otimes \phi_0 )\circ \rho_r ,\quad \psi = \rho_r^{-1} \circ (\mathrm{id} \otimes \psi_0 ).
\end{equation*}
where $ r=s+N$. Then for $S_\alpha z^k S_\beta^* \in \omega (s)$, we can show that
\begin{equation*}
\| \psi \circ \phi \circ \Phi^l (S_\alpha z^kS_\beta^* )- \Phi^l (S_\alpha z^kS_\beta^* )\|  <\delta
\end{equation*}
for $0 \le l \le N-1$. Hence,
\begin{eqnarray*}
\lefteqn{\limsup_{N \longrightarrow \infty} \frac{1}{N}\log rcp (\omega (s) \cup \Phi (\omega (s) ) \cup \cdots \cup \Phi^{N-1}(\omega (s)) ;\delta )}\hspace{1cm}\\
&\le& \limsup_{N \longrightarrow \infty} \frac{1}{N} \log ( Rn^r)
=\limsup_{N \longrightarrow \infty} \frac{1}{N}(\log R + (s+N) \log n)
=\log n
\end{eqnarray*}
and, using Theorem \ref{thm:KSthm}, we have finished the proof of $ht(\Phi ) \le \log n$. \\
On the other hand, we shall show that $\log n \le h_\varphi (\Phi )$. Using the gauge action of $\mathcal{O}_{(1,n)}(\mathbb{T}) $, we can take the conditional expectation onto the Bunce-Deddens algebra $\mathcal{B}_{(1,n)}$. Moreover, we consider $\mathbb{T}$-action on $\mathcal{B}_{(1,n)}$ defined as follow. First, we construct $\mathbb{T}$-action on $M_k(C(\mathbb{T} ))$ which are building blocks of $\mathcal{B}_{(1,n)}$. For $t \in \mathbb{T}$, let us define $\gamma^{(k)}_t : M_k(C(\mathbb{T} )) \longrightarrow M_k(C(\mathbb{T}))$ by
\begin{equation*}
\gamma^{(k)}_t(f)(z)=U_t^{(k)}f(t^kz)U_t^{(k)*} \quad (f \in M_k(C(\mathbb{T})) )
\end{equation*}
where $U_t^{(k)}$ is the unitary of $M_k(\mathbb{C} )$ defined by $U_t^{(k)}=\mathrm{diag}(1,t,\cdots ,t^{(k-1)} )$. Then these actions are compatible for the inductive limit system of $\mathcal{B}_{(1,n)}$, so we can construct the action arising from $\gamma^{(n^k)}$'s; we shall denote it by $\gamma$. Then we can check that the fixed point algebra $\mathcal{B}_{(1,n)}^\gamma$ is the continuous functions $C(K_n)$, where $K_n$ is the Cantor set, which is the maximal abelian algebra of Cuntz algebra. Hence, we obtain a conditional expectation onto $C(K_n)$ from $\mathcal{B}_{(1,n)}$ (and also from $\mathcal{O}_{(1,n)}(\mathbb{T}) $). Moreover, $\Phi|_{C(K_n)}$ is the canonical shift on $K_n$ and $\varphi|_{C(K_n)}$ is the canonical trace, so we obtain $h_{\varphi|_{C(K_n)}}(\Phi|_{C(K_n))}) = \log n$. Hence,
\begin{equation*}
\log n =h_{\varphi|_{C(K_n)}}(\Phi|_{C(K_n))}) \le h_{\varphi}(\Phi) \le ht(\Phi ) \le \log n.
\end{equation*}
Consequently $h_\varphi (\Phi )=ht(\Phi )=\log n$.

\end{proof}
\section{Dual action of the gauge action and $K$-theory}
\label{sec:Dual}

In \cite{Matsumoto1}, Matsumoto investigated the dual action of the gauge action on $C^*$-algebras associated with a subshift on the level of the $K$-groups to study dimension groups for the subshift. We follow his argument for $\mathcal{O}_{(m,n)}(\mathbb{T})$. Here, we compute the behavior of the dual action on $K$-groups.

Let $\alpha :\mathbb{T} \longrightarrow \mathcal{O}_{(m,n)}(\mathbb{T}) $ be the canonical gauge action and consider the crossed product $\mathcal{O}_{(m,n)}(\mathbb{T})  \rtimes_\alpha \mathbb{T}$, which is the universal $C^*$-algebra generated by the *-algebra $L^1(\mathbb{T},\mathcal{O}_{(m,n)}(\mathbb{T}) )$ whose multiplication and involution are defined as follows:
\begin{equation*}
f*g(t)=\int_\mathbb{T} f(s)\alpha_s(g(s^{-1}t))ds ,\quad f^*(t)=\alpha_t(f(t^{-1})^*).
\end{equation*}
for $f,g \in L^1(\mathbb{T} ,\mathcal{O}_{(m,n)}(\mathbb{T}) ),t \in \mathbb{T}$. Let $\widehat{\alpha}$ be the \textit{dual action} of $\alpha$ which is defined at the level of functions by $\widehat{\alpha}(f)(t)=tf(t)$.
The crossed product $\mathcal{O}_{(m,n)}(\mathbb{T})  \rtimes_\alpha \mathbb{T} \rtimes_{\widehat{\alpha}} \mathbb{Z}$ is stably isomorphic to $\mathcal{O}_{(m,n)}(\mathbb{T}) $.
 Let $p_0:\mathbb{T} \longrightarrow \mathcal{O}_{(m,n)}(\mathbb{T}) $ be the constant function whose value everywhere is the unit of $\mathcal{O}_{(m,n)}(\mathbb{T}) $.
  By \cite{Rosenberg}, the fixed point algebra $\mathcal{O}_{(m,n)}(\mathbb{T}) ^\alpha$ is isomorphic to the algebra $p_0(\mathcal{O}_{(m,n)}(\mathbb{T})  \rtimes_\alpha \mathbb{T})p_0$. The isomorphism between them is given by the correspondence $j: \mathcal{O}_{(m,n)}(\mathbb{T}) ^\alpha \ni x \longmapsto \widehat{x} \in L^1(\mathbb{T} ,\mathcal{O}_{(m,n)}(\mathbb{T})  ) \subset \mathcal{O}_{(m,n)}(\mathbb{T})  \rtimes_\alpha \mathbb{T}$ where the function $\widehat{x}$ is defined by $\widehat{x}(t)=x$ for $t \in \mathbb{T}$.

\begin{lem}\normalfont\slshape
The projection $p_0$ is full in $\mathcal{O}_{(m,n)}(\mathbb{T})  \rtimes_\alpha \mathbb{T}$.
\end{lem}

\begin{proof}
The proof of this lemma is the same as that of Lemma 4.1 of \cite{Matsumoto1}, but we give it for completeness. Suppose that there exists a nondegenerate representation $\pi$ of $\mathcal{O}_{(m,n)}(\mathbb{T}) \rtimes_\alpha \mathbb{T}$ such that $\pi (p_0)=0$. For any $x \in \mathcal{O}_{(m,n)}(\mathbb{T})$,
\begin{equation*}
\widehat{x}*p_0(t)=\int_\mathbb{T} \widehat{x}(s)\alpha_s (p_0(s^{-1}t))ds
=x.
\end{equation*}
Hence, $\widehat{x}*p_0 =\widehat{x}$. This implies that $\widehat{x} \in \ker \pi$.
For $x \in \mathcal{O}_{(m,n)}(\mathbb{T}) ,| \mu |=k \in \mathbb{N}$,
\begin{eqnarray*}
\widehat{xS_\mu}*\widehat{S_\mu}^*(t)
=\int_\mathbb{T} \widehat{xS_\mu}(s)\alpha_s( \widehat{S_\mu}^* (s^{-1}t))ds
=xS_\mu  \int_\mathbb{T} \alpha_s(\alpha_{s^{-1}t}(S_\mu^*)))ds
=xS_\mu \alpha_t(S_\mu^*)=t^{-k}xS_\mu S_\mu^*
\end{eqnarray*}
and we take the summation for the words of length $k$, 
$
\Bigl( \sum_{|\mu|=k} \widehat{xS_\mu}*\widehat{S_\mu}^* \Bigr)(t)
=t^{-k}x.
$
We can also show that
$
\widehat{xS_\mu^*}*\widehat{S_\mu^*}^*(t)=t^k x
$
 for $|\mu |=k$. 
Hence any $\mathcal{O}_{(m,n)}(\mathbb{T})$-valued function of the form 
$
\mathbb{T} \ni t \longmapsto t^k x
$ 
is contained in the ideal $\ker (\pi )$. This implies that $p_0$ is a full projection in $\mathcal{O}_{(m,n)}(\mathbb{T}) \rtimes_\alpha \mathbb{T}$.
\end{proof}

Since $p_0$ is a full projection in $\mathcal{O}_{(m,n)}(\mathbb{T})  \rtimes_\alpha \mathbb{T}$, there is an isometry $v \in M((\mathcal{O}_{(m,n)}(\mathbb{T})  \rtimes_\alpha \mathbb{T} ) \otimes \mathbb{K} )$ such that $v^*v=1 \otimes 1, v v^*=p_0 \otimes 1$ and 
\begin{equation*}
\mathrm{Ad}(v^*) : p_0 (\mathcal{O}_{(m,n)}(\mathbb{T})  \rtimes_\alpha \mathbb{T}) p_0 \otimes \mathbb{K} \longrightarrow (\mathcal{O}_{(m,n)}(\mathbb{T})  \rtimes_\alpha \mathbb{T}) \otimes \mathbb{K}
\end{equation*}
induce an isomorphism. We shall show that we can treat $\mathrm{Ad}(v^*)$ as an inclusion map $\iota :p_0 (\mathcal{O}_{(m,n)}(\mathbb{T})  \rtimes_\alpha \mathbb{T}) p_0 \otimes \mathbb{K} \longrightarrow (\mathcal{O}_{(m,n)}(\mathbb{T})  \rtimes_\alpha \mathbb{T}) \otimes \mathbb{K}$ in $K$-groups.

\begin{lem}\normalfont\slshape
$K_i(\iota )=K_i(\mathrm{Ad}(v^*) )$ for $i=0,1$.
\end{lem}

\begin{proof}
From Proposition 12.2.2 of \cite{Blackadar}, we can take a continuous path of isometries $(w_t)_{t \in (0,1]}$ in the multiplier algebra $M((\mathcal{O}_{(m,n)}(\mathbb{T})  \rtimes_\alpha \mathbb{T} ) \otimes \mathbb{K} )$ such that $w_tw_t^* \longrightarrow 0 \ (t \rightarrow 0)$ strictly. Put $v_t=w_tvw_t^*+(1-w_tw_t^*)$ for $t\in (0,1]$ and $v_0=1$. Then $( \mathrm{Ad} v_t^*(x) )_{t \in [0,1]}$ for $x \in p_0 (\mathcal{O}_{(m,n)}(\mathbb{T})  \rtimes_\alpha \mathbb{T}) p_0 \otimes \mathbb{K}$ is the norm continuous path i.e., $\mathrm{Ad}(v^*)$ and $\iota$ are homotopy equivalent. Hence, the above path implies that $K_i(\iota )=K_i(\mathrm{Ad}(v^*) )$ for $i=0,1$.
\end{proof}

First we consider the $K_0$-group. The group $K_0(\mathcal{O}_{(m,n)}(\mathbb{T}) ^\alpha )$ is isomorphic to $\mathbb{Z}[1/n]$ and for any $k \in \mathbb{N}$, $[S_1^kS_1^{k*}]_0$ corresponds to $1/n^k$. Note that $K_0(\mathcal{O}_{(m,n)}(\mathbb{T}) ^\alpha ) \cong K_0(\mathcal{O}_{(m,n)}(\mathbb{T})  \rtimes_\alpha \mathbb{T} )$ by the induced map $\psi_0=K_0(\iota \circ j) : [q]_0 \longmapsto [\widehat{q}]_0$. Let us define a map $\beta_0 : K_0(\mathcal{O}_{(m,n)}(\mathbb{T}) ^\alpha ) \longrightarrow K_0(\mathcal{O}_{(m,n)}(\mathbb{T}) ^\alpha )$ by $\beta_0=\psi_0^{-1} \circ K_0(\widehat{\alpha }) \circ \psi_0$:
\begin{equation*}
\begin{CD}
K_0(\mathcal{O}_{(m,n)}(\mathbb{T})  \rtimes_\alpha \mathbb{T} ) @ >K_0(\widehat{\alpha })>> K_0(\mathcal{O}_{(m,n)}(\mathbb{T})  \rtimes_\alpha \mathbb{T} ) \\
 @A\psi_0 A\cong A @ A\psi_0 A \cong A \\
K_0(\mathcal{O}_{(m,n)}(\mathbb{T}) ^\alpha ) @ >\beta_0 >> K_0(\mathcal{O}_{(m,n)}(\mathbb{T}) ^\alpha ).
\end{CD}
\end{equation*}

\begin{lem}\normalfont\slshape
\label{lem:beta0}
For any projection $q \in \mathcal{B}_{(m,n)}=\mathcal{O}_{(m,n)}(\mathbb{T}) ^\alpha,\ \beta_0[q]_0=[S_1qS_1^*]_0$.
\end{lem}

\begin{proof}
This proof is the same as Lemma 4.5 of \cite{Matsumoto1}, but for convenience, we repeat it.
It is enough to show that $K_0(\widehat{\alpha} )[\widehat{q}]_0=[\widehat{S_1qS_1^*}]_0$ in $K_0(\mathcal{O}_{(m,n)}(\mathbb{T})  \rtimes_\alpha \mathbb{T} )$. Since $q \in \mathcal{O}_{(m,n)}(\mathbb{T}) ^\alpha$, we have
\begin{eqnarray*}
(\widehat{S_1}*\widehat{q})(t)&=&\int_\mathbb{T} \widehat{S_1}(s)\alpha_s(\widehat{q}(s^{-1}t))ds=S_1 \int_\mathbb{T} \alpha_s(q)ds=S_1q\\
(\widehat{S_1}*\widehat{q}*\widehat{S_1}^*)(t)&=&
\int_\mathbb{T} (\widehat{S_1}*\widehat{q})(s)\alpha_s(\widehat{S_1}^*(s^{-1}t))ds=S_1q \alpha_t(S_1^*)=t^{-1}S_1qS_1^*.
\end{eqnarray*}
Hence, $\widehat{\alpha} (\widehat{S_1}*\widehat{q}*\widehat{S_1}^*)(t)=S_1qS_1^*$ and this implies that $\widehat{\alpha} (\widehat{S_1}*\widehat{q}*\widehat{S_1}^*)=\widehat{S_1qS_1^*}$. We can easily check that $\widehat{S_1}^* * \widehat{S_1}=p_0$. Put $W=\widehat{S_1}*\widehat{q}$; then $W^* *W=\widehat{q} * \widehat{S_1}^* * \widehat{S_1} * \widehat{q}=\widehat{q} * p_0 * \widehat{q}=\widehat{q},\ W * W^*=\widehat{S_1} * \widehat{q} * \widehat{S_1}^*$. This implies that 
\begin{equation*}
K_0( \widehat{\alpha} )[\widehat{q}]_0=K_0(\widehat{\alpha})[\widehat{S_1} * \widehat{q} * \widehat{S_1}^*]_0 = [\widehat{S_1qS_1^*}]_0 \mbox{ \quad in  $K_0(\mathcal{O}_{(m,n)}(\mathbb{T})  \rtimes \mathbb{T})$ }.
\end{equation*}
Hence, the proof is complete.
\end{proof}

Next, we consider the $K_1$-group. We remark that $K_1(\mathcal{O}_{(m,n)}(\mathbb{T}) ^\alpha )=\mathbb{Z}[1/m]$ and $[S_1^kzS_1^{*k}+1-S_1^kS_1^{*k}]_1$ corresponds to $1/m^k$. The map $\psi_1:=K_1(\iota \circ j) : K_1(\mathcal{O}_{(m,n)}(\mathbb{T}) ^\alpha )\longrightarrow K_1(\mathcal{O}_{(m,n)}(\mathbb{T})  \rtimes_\alpha \mathbb{T} )$ is determined by $\psi_1 [S_1^kzS_1^{k*}+1-S_1^kS_1^{k*}]_1=[(S_1^kzS_1^{k*}+1-S_1^kS_1^{k*})\widehat{\ }+1-p_0]_1$, where $1$ is the unit of the unitization $C^*$-algebra $(\mathcal{O}_{(m,n)}(\mathbb{T})  \rtimes_\alpha \mathbb{T})^\dagger$. Let us define $\beta_1 :K_1(\mathcal{O}_{(m,n)}(\mathbb{T}) ^\alpha ) \longrightarrow K_1(\mathcal{O}_{(m,n)}(\mathbb{T}) ^\alpha )$ by $\beta_1=\psi_1^{-1} \circ K_1(\widehat{\alpha} ) \circ \psi_1 ;$
\begin{equation*}
\begin{CD}
K_1(\mathcal{O}_{(m,n)}(\mathbb{T})  \rtimes_\alpha \mathbb{T} ) @ >K_1(\widehat{\alpha })>> K_1(\mathcal{O}_{(m,n)}(\mathbb{T})  \rtimes_\alpha \mathbb{T} ) \\
 @A\psi_1 A\cong A @ A\psi_1 A \cong A \\
K_1(\mathcal{O}_{(m,n)}(\mathbb{T}) ^\alpha ) @ >\beta_1 >> K_1(\mathcal{O}_{(m,n)}(\mathbb{T}) ^\alpha ).
\end{CD}
\end{equation*}
We recall the following lemma (Lemma 1.2 of \cite{Cuntz2}).

\begin{lem}\normalfont\slshape
\label{lem:Cuntzlemma}
Let $B$ be a $C^*$-algebra. Then for any partial isometry $s \in B^\dagger$ and unitary $u \in s^*s B^\dagger s^*s$, 
\begin{equation*}
[u+1-s^*s]_1=[sus^*+1-ss^*]_1 \mbox{ in $K_1(B)$}.
\end{equation*}
\end{lem}

Let us check that $\beta_1$ is the $1/m$-times map at the level of the $K_1$-group. It is enough to calculate for $[z]_1$.

\begin{lem}\normalfont\slshape
\label{lem:beta1}
$\beta_1[z]_1=\frac{1}{m}[z]_1 \mbox{ in } K_1(\mathcal{O}_{(m,n)}(\mathbb{T})^\alpha )=\mathbb{Z}[1/m]$.
\end{lem}

\begin{proof}
We compute $K_1(\widehat{\alpha} )[\widehat{z}+1-p_0]_1$ in $K_1(\mathcal{O}_{(m,n)}(\mathbb{T})  \rtimes_\alpha \mathbb{T} )$ instead of $\beta_1[z]_1$. From Lemma \ref{lem:Cuntzlemma}, $[\widehat{z}+(1-p_0)]_1=[\widehat{S_1} * \widehat{z} * \widehat{S_1}^* + 1-\widehat{S_1} * \widehat{S_1}^*]_1$. By a similar calculation to the one in Proposition \ref{lem:beta0},
\begin{eqnarray*}
\widehat{\alpha}^\dagger (\widehat{S_1} * \widehat{z} * \widehat{S_1}^* + 1-\widehat{S_1} * \widehat{S_1}^*)
&=&\widehat{S_1zS_1^*}+1-\widehat{S_1S_1^*}=\widehat{S_1zS_1^*}+p_0-\widehat{S_1S_1^*}+(1-p_0)\\
&=&(S_1zS_1^*+1-S_1S_1^*)\widehat{\ } +(1-p_0).
\end{eqnarray*}
This implies that
\begin{equation*}
\beta_1[z]_1=[S_1zS_1^*+(1-S_1S_1^*)]_1 \mbox{ \quad in $K_1(\mathcal{O}_{(m,n)}(\mathbb{T})  ^\alpha )$}.
\end{equation*}
Note that since $S_1 \notin \mathcal{O}_{(m,n)}(\mathbb{T} )^\alpha$, we cannot apply Lemma \ref{lem:Cuntzlemma}. The element $S_1zS_1^*+(1-S_1S_1^*)$ is $\mathrm{diag}(z,1_{n-1})$ in $\mathcal{B}_{(m,n)}$ and 
\[
\left[
\left(
\begin{array}{cccc|c}
z & 0 \\
0 & 1_{n-1}
\end{array}
\right)
\right]_1
=\frac{1}{m}
\left[
\left(
\begin{array}{cccc|c}
0 & z^m \\
1_{n-1} & 0
\end{array}
\right)
\right]_1
=\frac{1}{m}[z]_1.
\]
Hence, $\beta_1[z]_1=\frac{1}{m}[z]_1$.
\end{proof}

We summarize these lemmas below.

\begin{thm}
\label{thm:induced}
Suppose that $m \ge 1, n \ge 1$, and $\mathrm{gcd}(m,n)=1$. Let $\alpha$ be the gauge action of $\mathcal{O}_{(m,n)}(\mathbb{T})$ and $\widehat{\alpha}$ be the dual action of $\alpha$. For $i=0,1$, let $\psi_i : K_i(\mathcal{O}_{(m,n)}(\mathbb{T})^\alpha) \longrightarrow K_i(\mathcal{O}_{(m,n)}(\mathbb{T}) \rtimes_\alpha \mathbb{T})$ be the group isomorphism defined as above. Let $K_i(\widehat{\alpha})$ be the induced map of the dual action $\widehat{\alpha}$ on $K_i(\mathcal{O}_{(m,n)}(\mathbb{T}) \rtimes_\alpha \mathbb{T})$ and define $\beta_i:=\psi_i \circ K_i(\widehat{\alpha}) \circ \psi_i^{-1}$.\\
Then $\beta_0$ is a $1/n$-times map on $K_0(\mathcal{O}_{(m,n)}(\mathbb{T})^\alpha) \cong \mathbb{Z}[1/n]$ and $\beta_1$ is a $1/m$-times map on $K_1(\mathcal{O}_{(m,n)}(\mathbb{T})^\alpha) \cong \mathbb{Z}[1/m]$.
\end{thm}

We shall give another proof of the computation of $K$-groups of $\mathcal{O}_{(m,n)}(\mathbb{T}) $ (see also Proposition \ref{thm:Katsura}).
When we apply the Pimsner-Voiculescu six-term exact sequence for $\mathcal{O}_{(m,n)}(\mathbb{T})  \otimes \mathbb{K} =\mathcal{O}_{(m,n)}(\mathbb{T})  \rtimes_\alpha \mathbb{T} \rtimes_{\widehat{\alpha}} \mathbb{Z}$, we obtain the following exact sequence:
\[
\begin{CD}
K_0(\mathcal{O}_{(m,n)}(\mathbb{T})  \rtimes_\alpha \mathbb{T}) @>\mathrm{id}-K_0(\widehat{\alpha}^{-1})>> K_0(\mathcal{O}_{(m,n)}(\mathbb{T})  \rtimes_\alpha \mathbb{T}) @>>> K_0(\mathcal{O}_{(m,n)}(\mathbb{T})  )\\
@AAA @. @VVV \\
K_1(\mathcal{O}_{(m,n)}(\mathbb{T}) ) @<<< K_1(\mathcal{O}_{(m,n)}(\mathbb{T})  \rtimes_\alpha \mathbb{T}) @<\mathrm{id}-K_1(\widehat{\alpha}^{-1})<< K_1(\mathcal{O}_{(m,n)}(\mathbb{T})  \rtimes_\alpha \mathbb{T})
\end{CD}
\]
Since $K_i(\mathcal{O}_{(m,n)}(\mathbb{T})  \rtimes_\alpha \mathbb{T}) \cong K_i(\mathcal{O}_{(m,n)}(\mathbb{T}) ^\alpha)\ (i=0,1)$ and from the above argument, we have
\[
\begin{CD}
K_0(\mathcal{O}_{(m,n)}(\mathbb{T}) ^\alpha) @>\mathrm{id}-\beta_0^{-1}>> K_0(\mathcal{O}_{(m,n)}(\mathbb{T}) ^\alpha) @>>> K_0(\mathcal{O}_{(m,n)}(\mathbb{T})  )\\
@AAA @. @VVV \\
K_1(\mathcal{O}_{(m,n)}(\mathbb{T}) ) @<<< K_1(\mathcal{O}_{(m,n)}(\mathbb{T}) ^\alpha) @<\mathrm{id}-\beta_1^{-1}<< K_1(\mathcal{O}_{(m,n)}(\mathbb{T}) ^\alpha)
\end{CD}
\]
From $K_0(\mathcal{O}_{(m,n)}(\mathbb{T})^\alpha )=K_0(\mathcal{B}_{(m,n)})=\mathbb{Z}[1/n],K_1(\mathcal{O}_{(m,n)}(\mathbb{T})^\alpha )=K_1(\mathcal{B}_{(m,n)})=\mathbb{Z}[1/m]$ and Theorem \ref{thm:induced}, $\beta_0$ and $\beta_1$ correspond to the multiplication of $1/n$ and $1/m$, respectively. From an easy calculation, we obtain another computation of $K$-groups of $\mathcal{O}_{(m,n)}(\mathbb{T})$.


\section*{Acknowledgement}
\addcontentsline{toc}{section}{Acknowledgement}

The author thanks Professor Yasuo Watatani for constant encouragement and a lot of advices. He is also grateful to Kengo Matsumoto, who kindly taught me about his work.

\end{document}